\documentclass[11pt]{article}

\usepackage[english]{babel}
\usepackage{hyperref}

\usepackage[utf8]{inputenc}
\usepackage{amsfonts,amssymb,amsmath,amsthm,latexsym,epsfig,amscd,bbm,stmaryrd,mathrsfs}
\usepackage{graphicx,color,tikz,courier,bbm,enumerate,psfrag,wasysym,type1cm,dsfont,geometry,marginnote,afterpage}

\numberwithin{equation}{section}

\theoremstyle{plain} 

\newtheorem{theorem}{Theorem}[section]

\newtheorem{lemma}[theorem]{Lemma}
\newtheorem{proposition}[theorem]{Proposition}

\newtheorem{assumption}[theorem]{Assumption}
\newtheorem{remark}[theorem]{Remark}

\newcommand{\ssup}[1] {{{\scriptscriptstyle{{#1}}}}} 
\newcommand{\EE}{\mathbb{E}}
\newcommand{\PP}{\mathbb{P}}
\newcommand{\Prob}{\mathrm{P}}
\newcommand{\Probgr}{\mathfrak{P}}

\newcommand{\N}{\mathbb{N}}
\newcommand{\Z}{\mathbb{Z}}
\newcommand{\R}{\mathbb{R}}
\newcommand{\cO}{\mathcal{O}}
\newcommand{\cT}{\mathcal{T}}
\newcommand{\ee}{\mathrm{e}}
\newcommand{\dd}{\mathrm{d}}
\newcommand{\cP}{\mathcal{P}}
\newcommand{\GW}{\mathcal{GW}}
\newcommand{\supp}{\mathrm{supp}}
\newcommand{\scrX}{\mathscr{X}}

\newcommand{\BB}{\mathrm{BB}}
\newcommand{\bb}{\mathrm{bb}}
\newcommand{\cU}{\mathcal{U}}
\newcommand{\cS}{\mathcal{S}}
\newcommand{\cL}{\mathcal{L}}
\newcommand{\cW}{\mathcal{W}}

\begin{document}

\title{The annealed parabolic Anderson model\\ 
on a regular tree}

\author{
\renewcommand{\thefootnote}{\arabic{footnote}}
F.\ den Hollander
\footnotemark[1]
\\
\renewcommand{\thefootnote}{\arabic{footnote}}
D. Wang
\footnotemark[2]
}

\footnotetext[1]{
Mathematical Institute, Leiden University, P.O.\ Box 9512, 2300 RA Leiden, The Netherlands\\
{\tt denholla@math.leidenuniv.nl}
}

\footnotetext[2]{
Mathematical Institute, Leiden University, P.O.\ Box 9512, 2300 RA Leiden, The Netherlands\\
{\tt d.wang@math.leidenuniv.nl}
}

\date{\today}

\maketitle

\begin{abstract}
We study the total mass of the solution to the parabolic Anderson model on a regular tree with an i.i.d.\ random potential whose marginal distribution is double-exponential. In earlier work we identified two terms in the asymptotic expansion for large time of the total mass under the \emph{quenched law}, i.e., conditional on the realisation of the random potential. In the present paper we do the same for the \emph{annealed law}, i.e., averaged over the random potential. It turns out that the annealed expansion differs from the quenched expansion. The derivation of the annealed expansion is based on a \emph{new approach} to control the local times of the random walk appearing in the Feynman-Kac formula for the total mass. In particular, we condition on the backbone to infinity of the random walk, truncate and periodise the infinite tree relative to the backbone to obtain a random walk on a finite subtree with a specific boundary condition, employ the large deviation principle for the empirical distribution of \emph{Markov renewal processes} on finite graphs, and afterwards let the truncation level tend to infinity to obtain an asymptotically sharp asymptotic expansion.

\vspace{0.5cm}



\medskip\noindent
{\bf Acknowledgment:}
The research in this paper was supported by the Netherlands Organisation for Scientific Research through NWO Gravitation Grant NETWORKS-024.002.003. 

\normalsize
\end{abstract}

\newpage

\small
\tableofcontents
\normalsize

\maketitle

\newpage

%%%%%%%%%%%%%%%%%%%%%%%%%%%%%%%%%%%%%%%%%%%%%%%%

\section{Introduction and main results}

Section~\ref{motivation} provides background and motivation, Section~\ref{definitions} lists notations, definitions and assumptions, Section~\ref{maintheorem} states the main theorems, while Section~\ref{discussion} places these theorems in their proper context.

%%%

\subsection{Background and motivation}
\label{motivation}

The \emph{parabolic Anderson model} (PAM) is the Cauchy problem
\begin{equation}
\label{PAM}
\partial_t u(x,t) = \Delta_\scrX u(x,t) + \xi(x) u(x,t) , \qquad t>0, \, x \in \scrX,
\end{equation}
where $t$ is time, $\scrX$ is an ambient space, $\Delta_\scrX$ is a Laplace operator acting on functions on $\scrX$, and $\xi$ is a random potential on $\scrX$. Most of the literature considers the setting where $\scrX$ is either $\Z^d$ or $\R^d$ with $d \geq 1$, starting with the foundational papers \cite{GM1990}, \cite{GM1998}, \cite{GdH1999} and further developed through a long series of follow-up papers (see the monograph \cite{K2016} and the survey paper \cite{A2016} for an overview). More recently, other choices for $\scrX$ have been considered as well:
\begin{itemize} 
\item[(I)] 
\emph{Deterministic graphs} (the complete graph \cite{FM1990}, the hypercube \cite{AGH2016}).
\vspace{-.2cm}
\item[(II)] 
\emph{Random graphs} (the Galton-Watson tree \cite{dHKdS2020}, \cite{dHW2021}, the configuration model \cite{dHKdS2020}). 
\end{itemize} 
Much remains open for the latter class.  

The main target for the PAM is a description of \emph{intermittency}: for large $t$ the solution $u(\cdot,t)$ of \eqref{PAM} concentrates on well-separated regions in $\scrX$, called \emph{intermittent islands}. Much of the literature focusses on a detailed description of the size, shape and location of these islands, and on the profiles of the potential $\xi(\cdot)$ and the solution $u(\cdot,t)$ on them. A special role is played by the case where $\xi$ is an i.i.d.\ random potential with a \emph{double-exponential} marginal distribution 
\begin{equation}
\label{e:DEexact}
\Prob(\xi(0) > u) = \ee^{-\ee^{u/\varrho}}, \qquad u \in \R,
\end{equation}
where $\varrho \in (0,\infty)$ is a parameter. This distribution turns out to be critical, in the sense that the intermittent islands neither grow nor shrink with time, and represents a class of its own. 

In the present paper we consider the case where $\mathscr{X}$ is an \emph{unrooted regular tree} $\cT$. Our focus will be on the asymptotics as $t\to\infty$ of the total mass
\begin{equation*}
U(t) = \sum_{x \in \cT} u(x,t).
\end{equation*}
In earlier work \cite{dHKdS2020}, \cite{dHW2021} we were concerned with the case where $\mathscr{X}$ is a \emph{rooted Galton-Watson tree} in the \emph{quenched} setting, i.e., almost surely with respect to the random tree and the random potential. This work was restricted to the case where the random potential is given by \eqref{e:DEexact} and the offspring distribution of the Galton-Watson tree has support in $\N\backslash\{1\}$ with a sufficiently thin tail. In the present paper our focus will be on the \emph{annealed} setting, i.e., averaged over the random potential. We derive two terms in the asymptotic expansion as $t\to\infty$ of the average total mass 
\begin{equation*}
\langle U(t) \rangle = \sum_{x \in \cT} \langle u(x,t) \rangle,
\end{equation*} 
where $\langle\cdot\rangle$ denotes expectation with respect to the law of the random potential. It turns out that the annealed expansion \emph{differs} from the quenched expansion, even though the same variational formula plays a central role in the two second terms. 

The derivation in the annealed setting forces us to follow a different route than in the quenched setting, based on various approximations of $\cT$ that are more delicate than the standard approximation of $\Z^d$ (see \cite[Chapter VIII]{dHLDP2000}). This is the reason why we consider regular trees rather than Galton-Watson trees, to which we hope to return later. A key tool in the analysis is the large deviation principle for the empirical distribution of \emph{Markov renewal processes} on finite graphs derived in \cite{MZ2016}, which is recalled in Appendix \ref{appA}.

%%%

\subsection{The PAM on a graph}
\label{definitions}

%%%

\subsubsection{Notations and definitions} 

Let $G = (V,E)$ be a \emph{simple connected undirected} graph, either finite or countably infinite, with a designated vertex $\cO$ called the root. Let $\Delta_G$ be the Laplacian on $G$, i.e., 
\begin{equation*}
(\Delta_G f)(x) = \sum_{ {y\in V:} \atop { \{x,y\} \in E} } [f(y) - f(x)], \qquad x \in V,\,f\colon\,V\to\R,
\end{equation*}
which acts along the edges of $G$. Let $\xi := (\xi(x))_{x \in V}$ be a random potential attached to the vertices of $G$, taking values in $\R$. Our object of interest is the non-negative solution of the Cauchy problem with localised initial condition,
\begin{equation}
\label{e:PAMdef}
\begin{array}{llll}
\partial_t u(x,t) &=& (\Delta_G u)(x,t) + \xi(x) u(x,t), &x \in V,\, t>0,\\
u(x,0) &=& \delta_\cO(x), &x \in V.
\end{array}
\end{equation}
The quantity $u(x,t)$ can be interpreted as the amount of mass at time $t$ at site $x$ when initially there is unit mass at $\cO$. The total mass at time $t$ is $U(t) = \sum_{x \in V} u(x,t)$. The total mass is given by the \emph{Feynman-Kac formula}
\begin{equation}
\label{e:FKoriginal}
U(t) = \EE_\cO \left(\ee^{\int_0^t \xi(X_s) \dd s}\right),
\end{equation}
where $X=(X_t)_{t \geq 0}$ is the continuous-time random walk on the vertices $V$ with jump rate $1$ along the edges $E$, and $\PP_\cO$ denotes the law of $X$ given $X_0=\cO$. Let $\langle \cdot \rangle$ denote expectation with respect to $\xi$. The quantity of interest in this paper is the average total mass at time $t$:
\begin{equation}
\label{e:moments} 
\langle U(t) \rangle = \left\langle \EE_\cO\left(\ee^{\int_0^t \xi(X_s) \dd s}\right)\right\rangle. 
\end{equation}

%%%

\subsubsection{Assumption on the potential}

Throughout the paper we assume that the random potential $\xi = (\xi(x))_{x \in V}$ consists of i.i.d.\ random variables with a marginal distribution whose cumulant generating function 
\begin{equation}
\label{e:Hdef}
H(u) = \log \left\langle \ee^{u\xi(\cO)} \right\rangle
\end{equation}  
satisfies the following:

\begin{assumption}{\bf [Asymptotic double-exponential potential]} 
\label{ass:pot}
$\mbox{}$\\
{\rm There exists a $\varrho \in (0,\infty)$ such that
\begin{equation}
\label{e:DE}
\lim_{u\to\infty} u H''(u) = \varrho. 
\end{equation}
}\hfill$\spadesuit$
\end{assumption}

\begin{remark}{\bf [Double-exponential potential]} 
{\rm A special case of \eqref{e:DE} is when $\xi(\cO)$ has the double-exponential distribution in \eqref{e:DEexact}, in which case
\begin{equation*}
H(u) = \log \Gamma(\varrho u + 1)
\end{equation*}
with $\Gamma$ the gamma function.}\hfill$\spadesuit$
\end{remark}

By Stirling's approximation, \eqref{e:DE} implies
\begin{equation}
\label{e:Hasymp}
H(u) = \varrho u \log(\varrho u) - \varrho u + o(u), \qquad u \to \infty.
\end{equation}

Assumption~\ref{ass:pot} is more than enough to guarantee existence and uniqueness of the non-negative solution to \eqref{e:PAMdef} on any discrete graph with at most exponential growth (as can be inferred from the proof in \cite{GM1990}, \cite{GM1998} for the case $G=\Z^d$). Since $\xi$ is assumed to be i.i.d., we have from \eqref{e:moments} that
\begin{equation}
\label{e:FKann}
\langle U(t) \rangle = \mathbb{E}_\mathcal{O}\left(\exp\bigg[\sum\limits_{x\in V} H(\ell_t(x))\bigg]\right),
\end{equation}
where
\begin{equation*}
\ell_t(x) = \int^t_0 \mathds{1}\{X_s =x \}\, \dd s, \qquad x \in V,\, t\geq 0,
\end{equation*}
is the local time of $X$ at vertex $x$ up to time $t$.

%%%

\subsubsection{Variational formula}

The following \emph{characteristic variational formula} is important for the description of the asymptotics of $\langle U(t)\rangle$. Denote by $\cP(V)$ the set of probability measures on $V$. For $p \in \cP(V)$, define
\begin{equation}
\label{e:defIJ}
I_E(p) = \sum_{\{x,y\} \in E} \left( \sqrt{p(x)} - \sqrt{p(y)}\,\right)^2,
\qquad J_V(p) = - \sum_{x \in V} p(x) \log p(x),
\end{equation}
and set
\begin{equation}
\label{e:defchiG}
\chi_G(\varrho) = \inf_{p \in \cP(V)} [I_E(p) + \varrho J_V(p)], \qquad \varrho \in (0,\infty).
\end{equation}
The first term in \eqref{e:defchiG} is the quadratic form associated with the Laplacian, which is the large deviation rate function for
the \emph{empirical distribution}  
\begin{equation}
\label{e:emp}
L_t = \frac{1}{t} \int_0^t \delta_{X_s}\,\dd s = \frac{1}{t} \sum_{x \in V} \ell_t(x) \delta_x \in \cP(V)
\end{equation}
(see e.g.\ \cite[Section IV]{dHLDP2000}). The second term in \eqref{e:defchiG} captures the second order asymptotics of $\sum_{x \in V} H(tp(x))$ as $t \to\infty$ via \eqref{e:Hasymp} (see e.g.\ \cite[Section VIII]{dHLDP2000}). 

%%%

\subsubsection{Reformulation}

The following lemma pulls the leading order term out of the expansion and shows that the second order term is controlled by the large deviation principle for the empirical distribution. 

\begin{lemma}{\bf [Key object for the expansion]}
\label{lem:trunc}
If $G=(V,E)$ is finite, then
\begin{equation*}
\langle U(t) \rangle = \ee^{H(t) + o(t)}\,\EE_\cO\left(\ee^{-\varrho t J_V(L_t)}\right),
\qquad t \to \infty.      
\end{equation*}
where $J_V$ is the functional in \eqref{e:defIJ} and $L_t$ is the empirical distribution in \eqref{e:emp}.  
\end{lemma}

\begin{proof}
Because $\sum_{x \in V} \ell_t(x) = t$, we can rewrite \eqref{e:FKann} as
\begin{equation*}
\begin{aligned}
\langle U(t) \rangle &= \EE_\cO\left(\exp\bigg[\sum\limits_{x\in V} H(\ell_t(x))\bigg]\right)\\
&= \ee^{H(t)}\, \EE_\cO\left(\exp\bigg\{t \sum\limits_{x\in V} 
\frac{1}{t}\left[H(\tfrac{\ell_t(x)}{t}t)-\tfrac{\ell_t(x)}{t}H(t)\right]\bigg\}\right).
\end{aligned}
\end{equation*}
Assumption~\ref{ass:pot} implies that $H$ has the following scaling property (see \cite{GdH1999}):
\begin{equation*}
\lim_{t\to \infty} \frac{1}{t} [H(ct) - cH(t)] = \varrho c \log c \quad \text{uniformly in } c \in [0,1].
\end{equation*}
Hence the claim follows.
\end{proof}

%%%%%%%%%%%%%%%%%%%%%%%%%%%%%%%%%%%%%%%%%%%%%
\begin{figure}[htbp]
\begin{center}
\includegraphics[width=.35\textwidth]{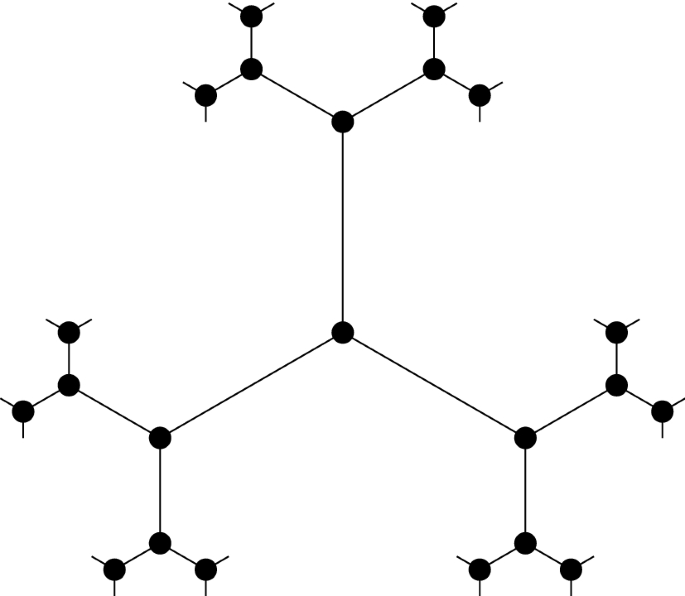}
\end{center}
\caption{\small An unrooted tree with degree $3$ (= offspring size $2$).}
\label{fig:tree}
\end{figure}
%%%%%%%%%%%%%%%%%%%%%%%%%%%%%%%%%%%%%%%%%%%%%%%%%

%%%

\subsection{The PAM on an unrooted regular tree: annealed total mass for large times and key variational formula}
\label{maintheorem}

In this section we specialise to the case where $G=\cT = (E,V)$, an unrooted regular tree of degree $d +1$ with $d \geq 2$ (see Fig.~\ref{fig:tree}). The main theorem of our paper is the following expansion.

\begin{theorem}{\bf [Growth rate of the total mass]}
\label{t:AMoments}
For any $d \geq 4$, subject to Assumption~{\rm \ref{ass:pot}},
\begin{equation}
\label{e:AMoments}
\frac{1}{t} \log \langle U(t) \rangle = \varrho \log(\varrho t) - \varrho  - \chi_\cT(\varrho) + o(1), \qquad t \to \infty,
\end{equation}
where $\chi_\cT(\varrho)$ is the variational formula in \eqref{e:defchiG} with $G=\cT$.   
\end{theorem}

\noindent
The proof of Theorem~\ref{t:AMoments} is given in Sections~\ref{proof:lb}--\ref{proof:ub} and makes use of technical computations collected in Appendices~\ref{appA}--\ref{appE}.

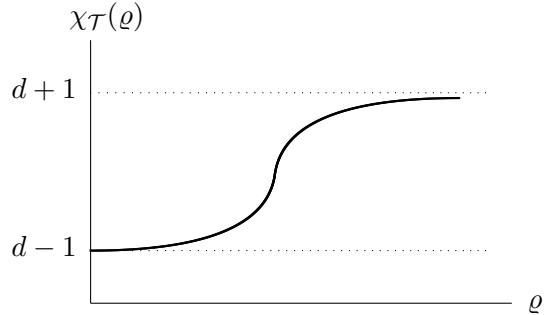
\begin{figure}[htbp]
\vspace{1cm}
\begin{center}
\setlength{\unitlength}{0.7cm}
\begin{picture}(8,5)(0,0)
%%%
\put(0,0){\line(1,0){8}}
\put(0,0){\line(0,1){5}}
%%%
{\thicklines
\qbezier(0,1)(3.3,1)(3.5,2.45)
\qbezier(3.5,2.45)(3.7,3.9)(7,3.9)
}
\qbezier[50](0,4)(4,4)(7.5,4)
\qbezier[50](0,1)(4,1)(7.5,1)
%%%
\put(8.3,-.1){$\varrho$}
\put(-.4,5.3){$\chi_\cT(\varrho)$}
\put(-1.5,0.9){$d-1$}
\put(-1.5,3.9){$d+1$}
%%%
\end{picture}
\end{center}
\caption{\small Qualitative plot of $\varrho \mapsto \chi_\cT(\varrho)$.}
\label{fig:chi}
\end{figure}

The main properties of the key quantity 
\begin{equation}
\label{varform}
\chi_\cT(\varrho) = \inf_{p \in \cP(V)} [I_E(p) + \varrho J_V(p)], \qquad \varrho \in (0,\infty),
\end{equation}
are collected in the following theorem (see Fig.~\ref{fig:chi}).

\begin{theorem}{\bf [Properties of the variational formula]}
\label{t:chivar}
For any $d \geq 2$ the following hold:\\
(a) The infimum in \eqref{varform} may be restricted to the set
\begin{equation}
\label{pdownset}
\cP_\cO^\downarrow(V) = \big\{p \in \cP(V)\colon\, \mathrm{argmax}\,p = \cO,\, 
p \text{ is non-increasing in the distance to } \cO\big\}.  
\end{equation}
(b) For every $\varrho \in (0,\infty)$, the infimum in \eqref{varform} restricted to $\cP_\cO^\downarrow(V)$  is attained, every minimiser $\Bar{p}$ is such that $\Bar{p}>0$ on $V$, and $\partial S_R = \sum_{\partial B_R(\cO)} \Bar{p}(x)$, $R\in\N_0$, satisfies
\begin{equation*}
\sum_{R \in \N_0} \partial S_R \log(R+1) \leq \frac{d+1}{\varrho}, 
\end{equation*}
where $B_R(\cO)$ is the ball of radius $R$ centred at $\cO$.\\
(c) The function $\varrho \mapsto \chi_\cT(\varrho)$ is strictly increasing and globally Lipschitz continuous on $(0,\infty)$, with 
\begin{equation*}
\lim_{\varrho \downarrow 0} \chi_\cT(\varrho) = d-1, \qquad \lim_{\varrho \to \infty} \chi_\cT(\varrho) = d+1.
\end{equation*} 
\end{theorem}

\noindent
The proof of Theorem~\ref{t:chivar} is given in Appendix~\ref{appC} (see Fig.~\ref{fig:chi}).

%%%%%%%%%%%%%%%%%%%%%%%%%%%%%%%%%%%%%%%%%%%%%
%%%%%%%%%%%%%%%%%%%%%%%%%%%%%%%%%%%%%%%%%%%%%%%%%

%%%

\subsection{Discussion}
\label{discussion}

{\bf 1.} 
Theorem~\ref{t:AMoments} identifies the scaling of the total mass up to and including terms that are exponential in $t$. The first two terms in the right-hand side of \eqref{e:AMoments} are the same as those of $\frac{1}{t} H(t)$ (recall \eqref{e:Hasymp}). The third term is a correction that comes from the cost for $X$ in the Feynman-Kac formula in \eqref{e:FKoriginal} to create an \emph{optimal local time profile} somewhere in $\cT$, which is captured by the minimiser(s) of the variational formula in \eqref{varform}.

\medskip\noindent
{\bf 2.} 
For the quenched model on a \emph{rooted Galton-Watson tree} $\GW$ we found in \cite{dHKdS2020}, \cite{dHW2021} that    
\begin{equation}
\label{e:AMomentsque}
\frac{1}{t} \log U(t) = \varrho \log\left(\frac{\varrho t \vartheta}{\log\log t}\right) 
- \varrho  - \chi(\varrho) +o(1), \qquad t \to \infty,
\qquad \Prob \times \Probgr\text{-a.s.},
\end{equation}
where $\Prob$ is the law of the potential, $\Probgr$ is the law of $\GW$, $\vartheta$ is the logarithm of the mean of the offspring distribution, and
\begin{equation}
\label{varinf}
\chi_\cT(\varrho) = \inf_{\cS \subset \GW} \chi_\cS(\varrho)
\end{equation} 
with $\chi_\cS(\varrho)$ given by \eqref{e:defchiG} and the infimum running over all subtrees of $\GW$. This result was shown to be valid as soon as the offspring distribution has support in $\N\backslash\{1\}$ (i.e., all degrees are at least 3) and has a sufficiently thin tail. The extra terms in \eqref{e:AMomentsque} come from the cost for $X$ in the Feynman-Kac formula in \eqref{e:FKoriginal} to travel in a time of order $o(t)$ to an optimal finite subtree with an optimal profile of the potential, referred to as \emph{intermittent islands}, located at a distance of order $\varrho t/\log\log t$ from $\cO$, and to subsequently spend most of its time on that subtree. In this cost the parameter $\vartheta$ appears, which is absent in \eqref{e:AMoments}. It was shown in \cite{dHKdS2020} that if $\varrho \geq 1/\log (d_{\min}+1)$, with $d_{\min}$ the minimum of the support of the offspring distribution, then the infimum in \eqref{varinf} is attained at the unrooted regular tree with degree $d_{\min}+1$, i.e., the \emph{minimal unrooted regular tree contained in} $\GW$, for which $\vartheta = \log d_{\min}$. Possibly the bound on $\varrho$ is redundant.

\medskip\noindent
{\bf 3.} In view of Lemma~\ref{lem:trunc} and the fact that Assumption~\ref{ass:pot} implies \eqref{e:Hasymp}, we see that the proof of Theorem~\ref{t:AMoments} amounts to showing that, on $\cT = (V,E)$, 
\begin{equation*}
\lim_{t\to\infty} \frac{1}{t} \log \EE_\cO\left(\ee^{-\varrho t J_V(L_t)}\right) = - \chi_\cT(\varrho).
\end{equation*}
We achieve this by deriving asymptotically matching upper and lower bounds. These bounds are obtained by truncating $\cT$ outside a ball of radius $R$, to obtain a finite tree $\cT_R$, deriving the $t\to\infty$ asymptotics for finite $R$, and letting $R\to\infty$ afterwards. For the lower bound we can use the standard truncation technique based on \emph{killing} $X$ when it exits $\cT_R$ and applying the large deviation principle for the empirical distribution of \emph{Markov processes} on finite graphs derived in \cite{DV75}. For the upper bound, however, we cannot use the standard truncation technique based on \emph{periodisation} of $X$ beyond radius $R$, because $\cT$ is an expander graph (see \cite[Chapter IV]{K2016} for a list of known techniques on $\Z^d$ and $\R^d$). Instead, we follow a route in which $\cT$ is approximated in successive stages by a version of $\cT_R$ with a \emph{specific boundary condition}, based on monitoring $X$ relative to its backbone to infinity. This route allows us to use the large deviation principle for the empirical distribution of \emph{Markov renewal processes} on finite graphs derived in \cite{MZ2016}, but we need the condition $d \geq 4$ to \emph{control} the specific boundary condition in the limit as $R \to\infty$ (see Remark~\ref{bdcontrol} for more details). The reason why the approximation of $\cT$ by finite subtrees is successful is precisely because in the parabolic Anderson model the \emph{total mass tends to concentrate on intermittent islands}.

\medskip\noindent
{\bf 4.} Theorem~\ref{t:chivar} shows that, modulo translations, the optimal strategy for $L_t$ as $t\to\infty$ is to be close to a minimiser of the variational formula in \eqref{varform} restricted to $\cP_\cO^\downarrow(V)$. Any minimiser is centred at $\cO$, strictly positive everywhere, non-increasing in the distance to $\cO$, and rapidly tending to zero. The following questions remain open:
\begin{itemize}
\item[(1)] 
Is the minimiser $\Bar{p}$ unique modulo translation? 
\item[(2)] 
Does $\Bar{p}(x)$ satisfy $\lim_{|x| \to \infty} |x|^{-1} \log \bar{p}(x) = -\infty$, with $|x|$ the distance between $x$ and $\cO$? 
\item[(3)] 
Is $\Bar{p}$ radially symmetric? 
\item[(4)] 
Is $\varrho \mapsto \chi_\cT(\varrho)$ analytic on $(0,\infty)$? 
\end{itemize}
We expect the answer to be yes for (1) and (2), and to be no for (3) and (4).    

%%%%%%%% SECTION 2 %%%%%%%%%%%%%%%%%%%%%%%%%%%%%%%

\section{Proof of the main theorem: lower bound}
\label{proof:lb}

In this section we prove the lower bound in Theorem~\ref{t:AMoments}, which is standard and straightforward. In Section~\ref{ss.killing} we obtain a lower bound in terms of a variational formula by \emph{killing} the random walk when it exits $\cT_R$. In Section~\ref{ss.varformlow} we derive the lower bound of the expansion by letting $R\to\infty$ in the variational formula. 

%%%

\subsection{Killing and lower variational formula}
\label{ss.killing}

Fix $R\in\mathbb{N}$. Let $\cT_R$ be the subtree of $\cT=(V,E)$ consisting of all the vertices that are within distance $R$ of the root $\cO$ and all the edges connecting them. Put $V_R=V_R(\cT_R)$ and $E_R = E(\cT_R)$. Let $\tau_R = \inf\{t \geq 0\colon\,X_t \notin V_R\}$ denote the first time that $X$ exits $\cT_R$. It follows from \eqref{e:FKann} that
\begin{equation*}
\langle U(t) \rangle \geq \EE_\cO\left(\exp\bigg[\sum\limits_{x\in V_R} 
H(\ell_t(x))\bigg]\mathds{1}\big\{\tau_R>t\big\}\right).  
\end{equation*}
Since $\cT_R$ is finite, Lemma~\ref{lem:trunc} gives 
\begin{equation*}
\langle U(t) \rangle \geq\ee^{H(t) + o(t)}\,\EE_\cO\left[\ee^{-\varrho t J_V(L_t)}
\mathds{1}\big\{\tau_R>t\big\}\right]   
\end{equation*}
with $J_V$ the functional defined in \eqref{e:defIJ}. As shown in \cite{G1977} (see also \cite{GM1998}), the family of sub-probability distributions $\PP_\cO(L_t \in \cdot\,, \tau_R>t)$, $t \geq 0$, satisfies the LDP on $\cP^R(V) = \{p \in \cP(V)\colon\,\mathrm{supp}(p) \subset V_R\}$ with rate function $I_E$, with $I_E$ the functional defined in \eqref{e:defIJ}. This is the \emph{standard} LDP for the empirical distribution of \emph{Markov processes}. Therefore, by Varadhan's Lemma, 
\begin{equation*}
\lim_{t\to\infty} \frac{1}{t} \log \EE_\cO\left[\ee^{-\varrho t J_V(L_t)}
\mathds{1}\big\{\tau_R>t\big\}\right]  = - \chi^-_R(\varrho)   
\end{equation*}
with
\begin{equation}
\label{chilb}
\chi^-_R(\varrho) = \inf_{p \in \cP^R(V)} [I_E(p) +\varrho J_V(p)],
\end{equation}
where we use that $p \mapsto J_V(p)$ is bounded and continuous (in the discrete topology) on $\cP^R(V)$. Note that   
\begin{equation*}
\lim_{t \to \infty} \frac{1}{t} \log \PP_\cO(\tau_R>t) = - \inf_{p\in \cP^R(V)} I_E(p) < 0,
\end{equation*}
which is non-zero because any $p \in \cP^R(V)$ is non-constant on $V$. The expression in \eqref{chilb} is the same as \eqref{e:defchiG} with $G=\cT$, except that $p$ is restricted to $V_R$.

%%%

\subsection{Limit of the lower variational formula}
\label{ss.varformlow}

Clearly, $R \mapsto \chi^-_R(\varrho)$ \text{ is non-increasing}. To complete the proof of the lower bound in Theorem~ \ref{t:AMoments}, it remains is to show the following.

\begin{lemma}
$\limsup_{R\to\infty} \chi^-_R(\varrho) \leq \chi_\cT(\varrho)$.    
\end{lemma}

\begin{proof}
Pick any $p \in \cP(V)$ such that $I_E(p)<\infty$ and $J_V(p)<\infty$. Let $p^{\ssup R}$ be the projection of $p$ onto $V_R$, i.e., 
\begin{equation*}
p^{\ssup R}(x) = \left\{\begin{array}{ll}
p(x), &x \in \mathrm{int}(V_R),\\
\sum_{y \geq x} p(y), &x \in \partial V_R,
\end{array}
\right.
\end{equation*}
where $y \geq x$ means that $y$ is an element of the progeny of $x$ in $\cT$. Since $p^{\ssup R} \in \cP^R(V)$, we have from \eqref{chilb} that $\chi^-_R(\varrho) \leq I_E(p^{\ssup R}) + \varrho J_V(p^{\ssup R})$. Trivially, $\lim_{R\to\infty} I_E(p^{\ssup R}) = I_E(p)$ and $\lim_{R\to\infty} J_V(p^{\ssup R}) = J_V(p)$, and so we have $\limsup_{R\to\infty} \chi^-_R(\varrho) \leq I_E(p) + \varrho J_V(p)$. Since this bound holds for arbitrary $p \in \cP(V)$, the claim follows from \eqref{chilb}.
\end{proof}

%%%%%%%% SECTION 3 %%%%%%%%%%%%%%%%%%%%%%%%%%%%%

\section{Proof of the main theorem: upper bound}
\label{proof:ub}

In this section we prove the upper bound in Theorem~\ref{t:AMoments}, which is more laborious and requires a more delicate approach than the standard periodisation argument used on $\Z^d$ . In Section~\ref{ss.projper} we obtain an upper bound in terms of a variational formula on a version of $\cT_R$ with a specific boundary condition. The argument comes in four steps, encapsulated in Lemmas~\ref{lem1}--\ref{lem6} below:
\begin{itemize} 
\item[(I)] 
\emph{Condition on the backbone} of $X$ (Section \ref{sss.1}). 
\item[(II)] 
\emph{Project} $X$ onto a concatenation of finite subtrees attached to this backbone that are \emph{rooted} versions of $\cT_R$ (Section \ref{sss.2}). 
\item[(III)] 
\emph{Periodise} the projected $X$ to obtain a \emph{Markov renewal process} on a single finite subtree and show that the periodisation can be chosen such that the local times at the vertices on the \emph{boundary} of the finite subtree are negligible (Section \ref{sss.3}).
\item[(IV)] 
Use the \emph{large deviation principle} for the empirical distribution of \emph{Markov renewal processes} derived in \cite{MZ2016} to obtain a variational formula on a single subtree (Section \ref{sss.4}).
\end{itemize}
In Section~\ref{ss.varformup} we derive the upper bound of the expansion by letting $R\to\infty$ in the variational formula. 

%%%

\subsection{Backbone, projection, periodisation and upper variational formula}
\label{ss.projper}

%%%

\subsubsection{Backbone}
\label{sss.1}

For $r \in \N_0$, let $\tau_r$ be the last time when $X$ visits $\partial B_r(\cO)$, the boundary of the ball of radius $r$ around $\cO$. Then the sequence $\BB = (X_{\tau_r})_{r \in \N_0}$ forms the backbone of $X$, running from $\cO$ to infinity.

\begin{lemma}{\bf [Condition on a backbone]}
\label{lem1}
For every backbone $\bb$ and every $t \geq 0$,
\begin{equation*}
\mathbb{E}_\mathcal{O}\left(\exp\bigg[\sum\limits_{x\in V(\cT)} H(\ell_t(x))\bigg]\right)
= \mathbb{E}_\mathcal{O}\left(\exp\bigg[\sum\limits_{x\in V(\cT)} H(\ell_t(x))\bigg] ~\Bigg|~ \BB = \bb\right).
\end{equation*}
\end{lemma}

\begin{proof}
By symmetry, the conditional expectation in the right-hand side does not depend on the choice of $\bb$. Indeed, permutations of the edges away from the root do not affect the law of $\sum_{x\in V(\cT)} H(\ell_t(x))$.  
\end{proof}

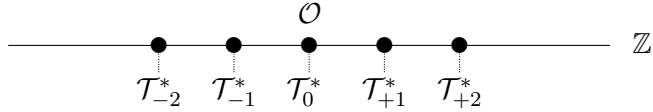
\begin{figure}[htbp]
\begin{center}
\setlength{\unitlength}{1cm}
\begin{picture}(8,1)(2,0)
%%%
\put(1,0){\line(1,0){8}}
%%%
\qbezier[10](3,0)(3,-.2)(3,-.35)
\qbezier[10](4,0)(4,-.2)(4,-.35)
\qbezier[10](5,0)(5,-.2)(5,-.35)
\qbezier[10](6,0)(6,-.2)(6,-.35)
\qbezier[10](7,0)(7,-.2)(7,-.35)
%%%
\put(2.7,-.7){$\cT^\ast_{-2}$}
\put(3.7,-.7){$\cT^\ast_{-1}$}
\put(4.7,-.7){$\cT^\ast_{0}$}
\put(5.7,-.7){$\cT^\ast_{+1}$}
\put(6.7,-.7){$\cT^\ast_{+2}$}
\put(3,0){\circle*{.2}}
\put(4,0){\circle*{.2}}
\put(5,0){\circle*{.2}}
\put(6,0){\circle*{.2}}
\put(7,0){\circle*{.2}}
\put(9.3,-.1){$\Z$}
\put(4.85,0.3){$\cO$}
%%%
\end{picture}
\vspace{1cm}
\end{center}
\caption{\small Redrawing of $\cT$ as $\cT^\Z$: a two-sided backbone $\Z$ with a sequence $\cT^*= (\cT^\ast_x)_{x \in \Z}$ of rooted trees hanging off. The upper index $*$ is used to indicate that the tree is rooted.}
\label{fig:redrawing}
\end{figure}

Turn the one-sided backbone into a two-sided backbone by adding a second backbone from $\cO$ to infinity. By symmetry, the choice of this second backbone is arbitrary, say $\bb'$.  Redraw $\cT$ by representing $\bb' \cup \bb$ as $\Z$ and representing the rest of $\cT$ as a sequence of \emph{rooted trees} $\cT^\ast = (\cT^\ast_x)_{x \in \Z}$ hanging off $\Z$  (see Fig.~\ref{fig:redrawing}). In $\cT^\ast_x$, the root sits at $x$ and has $d-1$ downward edges, while all lower vertices have $d$ downward edges.

%%%%%%%%%%%%%%%%%%%%%%%%%%%%%%%%%%%%%%%%%%%%%
%%%%%%%%%%%%%%%%%%%%%%%%%%%%%%%%%%%%%%%%%%%%%%%%%

Let $X^\Z=(X^\Z_t)_{t \geq 0}$ be the random walk on $\cT^\Z$ and $(\ell^\Z_t(x))_{x \in \cT^\Z}$ the local times of $X^\Z$ at time $t$. 

\begin{lemma}{\bf [Representation of $\cT$ as a backbone with rooted trees]}
For every $\bb$ and $t \geq 0$,
\begin{equation*}
\mathbb{E}_\mathcal{O}\left(\exp\bigg[\sum\limits_{x\in V(\cT)} H(\ell_t(x))\bigg] ~\Bigg|~ \BB = \bb\right)
= \mathbb{E}_\mathcal{O}\left(\exp\bigg[\sum\limits_{x\in V(\cT^{\,\Z})} H(\ell^\Z_t(x))\bigg] 
~\Bigg|~ X^\Z_\infty = + \infty\right).  
\end{equation*}
\end{lemma}

\begin{proof}
Simply redraw $\cT$ as $\cT^\Z$.
\end{proof}

\noindent
Note that $X^\Z$ is a Markov process whose sojourn times have distribution $\mathrm{EXP}(d+1)$ and whose steps are drawn uniformly at random from the $d+1$ edges that are incident to each vertex. 

%%%

\subsubsection{Projection}
\label{sss.2}

\begin{figure}[htbp]
\begin{center}
\setlength{\unitlength}{0.8cm}
\begin{picture}(7,7)(-1,-1)
%%%
\qbezier[50](0,0)(2.5,0)(5,0)
\qbezier[50](0,0)(1.25,2.5)(2.5,5)
\qbezier[50](5,0)(3.75,2.5)(2.5,5)
{\thicklines
\qbezier(2.5,5)(2.25,5.5)(2.0,6)
\qbezier(5,0)(5.25,-.5)(5.5,-1)
\qbezier(0,0)(0,-.5)(0,-1)
\qbezier(1,0)(1,-.5)(1,-1)
\qbezier(4,0)(4,-.5)(4,-1)
}
%%%
\put(0,0){\circle*{.2}}
\put(2.5,5){\circle*{.2}}
\put(5,0){\circle*{.2}}
\put(1,0){\circle*{.2}}
\put(4,0){\circle*{.2}}
%%%
\put(-0.17,-1.25){$\Box$}
\put(0.83,-1.25){$\Box$}
\put(3.83,-1.25){$\Box$}
%%%
\put(2.2,2){$\cT^*_R$}
%%%
\end{picture}
\end{center}
\caption{\small A unit in $\cU_R$. Inside is a rooted tree $\cT^*_R$ of depth $R-1$, of which only the root and the leaves are drawn. Hanging off the leaves at depth $R-1$ from the root are tadpoles, except for the right-most bottom vertex, which has a downward edge that connects to the root of the next unit. The vertices marked by a bullet form the boundary of $\cU_R$, the vertices marked by a square box form the tadpoles of $\cU_R$.}
\label{fig:unit}
\end{figure}
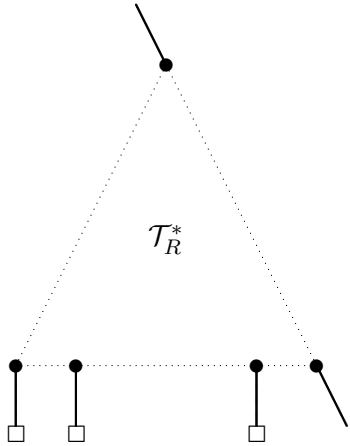

For $R \in \N\backslash\{1\}$, cut $\Z$ into slices of length $R$, i.e.,
\begin{equation*}
\Z = \cup_{k\in\Z} (z + (kR+I)), \qquad I=\{0,1,\ldots,R-1\},
\end{equation*} 
where $z$ is to be chosen later. Apply the following two maps to $\cT^\Z$ (in the order presented):
\begin{itemize}
\item[(i)]
For each $k \in \Z$, fold $\cT^\ast_{z+(kR+(R-1))}$ onto $\cT^\ast_{z+(k+1)R}$ by folding the $d-1$ edges downwards from the root on top of the edge in $\Z$ connecting $z+(kR+(R-1))$ and $z+(k+1)R$, and putting the $d$ infinite rooted trees hanging off each of these $d-1$ edges on top of the rooted tree $\cT^*_{z+(k+1)R}$ hanging off $z+(k+1)R$. Note that each of the $d$ infinite rooted trees is a copy of $\cT^*_{z+(k+1)R}$.
\item[(ii)]
For each $k \in \Z$ and $m \in \{0,1,\ldots,R-2\}$, cut off all the infinite subtrees trees in $\cT^\ast_{z+(kR+m)}$ whose roots are at depth $(R-1)-m$. Note that the total number of leaves after the cutting equals 
\begin{equation*}
(d-1) \sum_{m=0}^{R-2} d^{(R-2)-m} = (d-1)d^{R-2}\,\frac{1-d^{-(R-1)}}{1-d^{-1}} = d^{R-1} - 1,
\end{equation*} 
which is the same as the total number of leaves of the rooted tree $\cT^*_R$ of depth $R-1$ (i.e., with $R$ generations) minus $1$ (a fact we will need below).  
\end{itemize}
By doing so we obtain a \emph{concatenation} of finite units 
\begin{equation*}
\cU_R=(\cU_R[k])_{k \in \Z}
\end{equation*}
that are rooted trees of depth $R-1$ (see Fig.~\ref{fig:unit}). Together with the two maps that turn $\cT^\Z$ into $\cU_R$, we apply two maps to $X^\Z$:
\begin{itemize}
\item[(i)]
All excursions of $X^\Z$ in the infinite subtrees that are \emph{folded to the right and on top} are projected accordingly.
\item[(ii)]
All excursions of $X^\Z$ in the infinite subtrees that are \emph{cut off} are replaced by a sojourn of $X^{\cU_R}$ in the \emph{tadpoles} that replace these subtrees (see Fig.~\ref{fig:unit})     
\end{itemize}
The resulting path, which we call $X^{\cU_R} = (X^{\cU_R}_t)_{t \geq 0}$, is a \emph{Markov renewal process} with the following properties:
\begin{itemize}
\item
The sojourn times in all the vertices that are not tadpoles have distribution $\mathrm{EXP}(d+1)$.  
\item
The sojourn times in all the tadpoles have distribution $\psi$, defined as the conditional distribution of the \emph{return time} $\tau$ of the random walk on the infinite rooted tree $\cT^*$ \emph{given} that $\tau<\infty$ (see \cite{LP2016} for a proper definition).
\item
The transitions into the tadpoles have probability $\frac{d}{d+1}$, the transitions out of the tadpoles have probability $1$ (because of the condition $X^\Z_\infty = + \infty$).
\item
The transitions from $z + (kR+(R-1))$ to $z+(k+1)R$ have probability $\frac{d}{d+1}$, while the reverse transitions have probability $\frac{1}{d+1}$.    
\end{itemize} 
Write $(\ell^{\,\cU_R}_t(x))_{x \in V_{\cU_R}}$ to denote the local times of $X^{\cU_R}$ at time $t$.   
 
\begin{lemma}{\bf [Projection onto a concatenation of finite subtrees]} 
\label{lem3} 
For every $R \in \N\backslash\{1\}$ and $t \geq 0$,
\begin{equation*}
\begin{aligned}
&\mathbb{E}_\mathcal{O}\left(\exp\bigg[\sum\limits_{x\in V(\cT^{\,\Z})} H(\ell^\Z_t(x))\bigg]
~\Bigg|~ X^\Z_\infty = + \infty\right)\\
&\qquad \leq \mathbb{E}_\mathcal{O}\left(\exp\bigg[\sum\limits_{x\in V(\cU_R)} H(\ell^{\,\cU_R}_t(x))\bigg] 
~\Bigg|~ X^{\cU_R}_\infty = + \infty\right).
\end{aligned}
\end{equation*}
\end{lemma}

\begin{proof}
The maps that are applied to turn $X^\Z$ into $X^{\cU_R}$ are such that local times are \emph{stacked on top of each other}. Since $H$ defined in \eqref{e:Hdef} is convex and $H(0)=0$, we have $H(\ell) + H(\ell') \leq H(\ell+\ell')$ for all $\ell,\ell' \in \N_0$, which implies the inequality. 
\end{proof}

%%%%%%%%%%%%%%%%%%%%%%%%%%%%%%%%%%%%%%%%%%%%%
%%%%%%%%%%%%%%%%%%%%%%%%%%%%%%%%%%%%%%%%%%%%%%%%%

%%%

\subsubsection{Periodisation}
\label{sss.3}

Our next observation is that the condition $\{X^{\cU_R}_\infty = + \infty\}$ is \emph{redundant}. 

\begin{lemma}{\bf [Condition redundant]}
For every $R \in \N\backslash\{1\}$ and $t \geq 0$,
\begin{equation*}
\mathbb{E}_\mathcal{O}\left(\exp\bigg[\sum\limits_{x\in V(\cU_R)} H(\ell^{\,\cU_R}_t(x))\bigg]\,
~\Bigg|~ X^{\cU_R}_\infty = + \infty\right)
=  \mathbb{E}_\mathcal{O}\left(\exp\bigg[\sum\limits_{x\in V(\cU_R)} H(\ell^{\,\cU_R}_t(x))\bigg]\,\right).
\end{equation*}
\end{lemma}

\begin{proof}
The event $\{X^{\cU_R}_\infty = + \infty\}$ has probability $1$ because on the edges connecting the units of $\cU_R$ (see Fig.~\ref{fig:unit}) there is a drift downwards. To see why, note that $\frac{1}{d+1} < \tfrac12 < \frac{d}{d+1}$ because $d \geq 2$, and use that a one-dimensional random walk with drift is transient to the right \cite{S1976}.
\end{proof}

\begin{figure}[htbp]
\begin{center}
\setlength{\unitlength}{0.8cm}
\begin{picture}(7,7)(-1,-2)
%%%
\qbezier[50](0,0)(2.5,0)(5,0)
\qbezier[50](0,0)(1.25,2.5)(2.5,5)
\qbezier[50](5,0)(3.75,2.5)(2.5,5)
{\thicklines
\qbezier(2.5,5)(6,2.5)(5,0)
\qbezier(0,0)(0,-.5)(0,-1)
\qbezier(1,0)(1,-.5)(1,-1)
\qbezier(4,0)(4,-.5)(4,-1)
}
%%%
\put(0,0){\circle*{.2}}
\put(2.5,5){\circle*{.2}}
\put(5,0){\circle*{.2}}
\put(1,0){\circle*{.2}}
\put(4,0){\circle*{.2}}
%%%
\put(-0.17,-1.25){$\Box$}
\put(0.83,-1.25){$\Box$}
\put(3.83,-1.25){$\Box$}
%%%
\put(2.2,2){$\cT^*_R$}
%%%
\end{picture}
\vspace{-.8cm}
\end{center}
\caption{\small A unit $\cW_R$ with the top vertex and the right-most bottom vertex connected by an edge.}
\label{fig:foldedunit}
\end{figure}
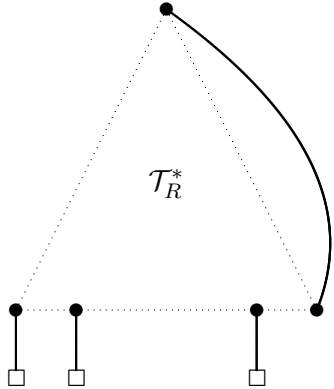

Since $\cU_R$ is periodic, we can \emph{fold} $X^{\cU_R}$ onto a single unit $\cW_R$, to obtain a Markov renewal process $X^{\cW_R}$ on $\cW_R$ (see Fig.~\ref{fig:foldedunit}) in which the transition from the top vertex to the right-most bottom vertex has probability $\frac{1}{d+1}$, while the reverse transition has probability $\frac{d}{d+1}$. Clearly, the sojourn time distributions are not affected by the folding and therefore remain as above. Write $(\ell^{\,\cW_R}_t(x))_{x \in V(\cW_R)}$ to denote the local times of $X^{\cW_R}$ at time $t$.   

\begin{lemma}{\bf [Periodisation to a single finite subtree]} 
For every $R \in \N\backslash\{1\}$ and $t \geq 0$,
\begin{equation*}
\mathbb{E}_\mathcal{O}\left(\exp\bigg[\sum\limits_{x\in V(\cU_R)} H(\ell^{\,\cU_R}_t(x))\bigg]\right)
\leq \mathbb{E}_\mathcal{O}\left(\exp\bigg[\sum\limits_{x\in V(\cW_R)} H(\ell^{\,\cW_R}_t(x))\bigg]\right).
\end{equation*}
\end{lemma}

\begin{proof}
The periodisation again stacks local time on top of each other.
\end{proof}

Before we proceed we make a \emph{crucial observation}, namely, we may still choose the shift $z \in \{0,1,\ldots,R-1\}$ of the cuts of the two-sided backbone $\Z$ (recall Fig.~\ref{fig:redrawing}). We will do so in such a way that the \emph{local time} up to time $t$ spent in the set $\partial_{\,\cU_R}$ defined by
\begin{equation}
\label{parRdef}
\begin{aligned}
\partial_{\,\cU_R} &= \text{all vertices at the top or at the bottom of a unit in $\cU_R$}\\
&= \text{all vertices marked by $\bullet$ in Fig.~\ref{fig:unit}}  
\end{aligned}
\end{equation}
is at most $t/R$. After the periodisation these vertices are mapped to the set  $\partial_{\,\cW_R}$ defined by
\begin{equation*}
\begin{aligned}
\partial_{\,\cW_R} &= \text{all vertices at the top or at the bottom of $\cW_R$}\\
&= \text{all vertices marked by $\bullet$ in Fig.~\ref{fig:foldedunit}}.  
\end{aligned}
\end{equation*}

\begin{lemma}{\bf [Control on the time spent at the boundary]}
\label{lem6}
For every $R \in \N\backslash\{1\}$ and $t \geq 0$,
\begin{equation*}
\begin{aligned}
&\mathbb{E}_\mathcal{O}\left(\exp\bigg[\sum\limits_{x\in V(\cU_R)} H(\ell^{\,\cU_R}_t(x))\bigg]\right)\\
&\qquad \leq \mathbb{E}_\mathcal{O}\left(\exp\bigg[\sum\limits_{x\in V(\cW_R)} H(\ell^{\,\cW_R}_t(x))\bigg]
\,1_{\big\{\frac{1}{t} \sum_{x \in \partial_{\,\cW_R}} \ell^{\cW_R}_t(x) \leq 1/R\big\}}\right).
\end{aligned}
\end{equation*}
\end{lemma}

\begin{proof}
For different $z$ the sets of vertices making up $\partial_R$ correspond to \emph{disjoint} sets of vertices in $\cT^\Z$ (see Fig.~\ref{fig:unit}). Since $\sum_{x \in \cT^\Z} \ell^\Z_t(x) = t$ for all $t \geq 0$, it follows that there exists a $z$ for which $\sum_{x \in \partial_R} \ell^\Z_t(x) \leq t/R$. Therefore the upper bound in Lemma~\ref{lem3} can be strengthened to the one that is claimed.
\end{proof}

%%%%%%%%%%%%%%%%%%%%%%%%%%%%%%%%%%%%%%%%%%%%%
%%%%%%%%%%%%%%%%%%%%%%%%%%%%%%%%%%%%%%%%%%%%%%%%% 

%%%

\subsubsection{Upper variational formula} 
\label{sss.4}

Lemmas~\ref{lem1}--\ref{lem6} provide us with an upper bound for the average total mass (recall (\eqref{e:FKann}) on the \emph{infinite} tree $\cT$ in terms of the same quantity on the \emph{finite} tree-like unit $\cW_R$ with a \emph{specific boundary condition}. Along the way we have paid a price: the sojourn times in the tadpoles are \emph{no longer} exponentially distributed, and the transition probabilities into and out of the tadpoles and between the top vertex and the right-most bottom vertex are \emph{biased}. We therefore need the large deviation principle for the empirical distribution of Markov renewal processes derived in \cite{MZ2016}, which we can now apply to the upper bound.

Since $\cW_R$ is finite, Lemma~\ref{lem:trunc} gives 
\begin{equation*}
\langle U(t) \rangle \leq \ee^{H(t) + o(t)}\,
\mathbb{E}_\mathcal{O}\left(\ee^{-\varrho J_{V(\cW_R)}(L^{\,\cW_R}_t)}\,
1_{\big\{L^{\cW_R}_t(\partial_{\,\cW_R}) \leq 1/R\big\}}\right)
\end{equation*}
with $J_V$ the functional defined in \eqref{e:defIJ}. The following lemma controls the expectation in the right-hand side.

\begin{lemma}{\bf [Scaling of the key expectation]}
\label{lem:VF}
For every $R \in \N\backslash\{1\}$,
\begin{equation*}
\lim_{t\to\infty} \frac{1}{t} \log \EE_\cO\left(\ee^{-\varrho t J_{V(\cW_R)}(L^{\cW_R}_t)}
\,1_{\big\{L^{\cW_R}_t(\partial_{\,\cW_R}) \leq 1/R\big\}}\right) = - \chi^+_R(\varrho),
\end{equation*}
where 
\begin{equation}
\label{VF}
\chi^+_R(\varrho) = \inf_{ {p \in \cP(V(\cW_R))\colon} \atop {p(\partial_{\cW_R}) \leq 1/R}}
\left\{I^\dagger_{E(\cW_R)}(p) + \varrho J_{V(\cW_R)}(p)\right\},
\end{equation}
with 
\begin{equation}
\label{RFup}
I^\dagger_{E(\cW_R)}(p) = \inf_{\beta \in (0,\infty)} \inf_{q \in \cP(V(\cW_R))} \big[\widehat{K}(\beta q) + \widetilde{K}(p \mid \beta q)\big],
\end{equation}
where
\begin{eqnarray}
\label{rf1}
\widehat{K}(\beta q) &=& \sup_{\widehat{q} \in \cP(V(\cW_R))} \sum_{x \in V(\cW_R)} 
\beta q(x) \log\left(\tfrac{\widehat{q}(x)}{\sum_{y \in V(\cW_R)} \pi_{x,y}\widehat{q}(y)}\right),\\
\label{rf2}    
\widetilde{K}(p \mid \beta q) &=& \sum_{x \in V(\cW_R)} \beta q(x)\,(\cL\lambda_x)\left(\tfrac{p(x)}{\beta q(x)}\right),
\end{eqnarray}
with
\begin{eqnarray}
\label{Legendre}
(\cL\lambda_x)(\alpha) &=&  \sup_{\theta \in \mathbb{R}} [\alpha\theta - \lambda_x(\theta)], \quad\quad \alpha \in [0,\infty),\\
\label{cumulant}
\lambda_x(\theta) &=& \log \int_0^\infty \ee^{\theta \tau}\psi_x(\dd\tau), \quad \theta \in \mathbb{R},
\end{eqnarray}
where $\psi_x=\psi$ when $x$ is a tadpole, $\psi_x = \mathrm{EXP}(d+1)$ when $x$ is not a tadpole, and $\pi_{x,y}$ is the transition kernel of the discrete-time Markov chain on $V(\cW_R)$ embedded in $X^{\cW_R}$.
\end{lemma}

\begin{proof}
Apply the large deviation principle derived in \cite{MZ2016}, which we recall in Proposition~\ref{prop:LDPrenewal} in Appendix~\ref{appA}. 
\end{proof}

The expression in \eqref{VF} is similar to \eqref{e:defchiG} with $G=\cW_R$, expect that the rate function $I_{E(\cW_R)}$ in \eqref{RFup} is more involved than the rate function $I_E$ in \eqref{e:defIJ}.  

%%%

\subsection{Limit of the upper variational formula}
\label{ss.varformup}

The prefactor $\ee^{H(t)+o(1)}$ in Lemma~\ref{lem:trunc} accounts for the terms $\varrho\log(\varrho t)-\varrho$ in the right-hand side of \eqref{e:AMoments} (recall \ref{e:Hasymp}). In view of Lemma~\ref{lem:VF}, in order to complete the proof of the upper bound in Theorem~\ref{t:AMoments} it suffices to prove the following lemma.

\begin{lemma}
For any $d \geq 4$, $\liminf_{R\to\infty} \chi^+_R(\varrho) \geq \chi_\cT(\varrho)$.
\end{lemma}

\begin{proof}
The proof is given in Appendix~\ref{appE} and relies on two steps:
\begin{itemize} 
\item
Show that, for $d \geq 4$, 
\begin{equation}
\label{IElberror}
I^\dagger_{E(\cW_R)}(p) \geq I^+_{E(\cW_R)}(p) + O(1/R)
\end{equation} 
with $I^+_{E(\cW_R)}$ a rate function similar to the \emph{standard rate function} $I_{E(\cW_R)}$ given by \eqref{e:defIJ}.
\item
Show that, $d \geq 2$, 
\begin{equation*}
\widehat{\chi}^{\,+}_R(\varrho) = \inf_{ {p \in \cP(V(\cW_R))\colon} 
\atop {p(\partial_{\,\cW_R}) \leq 1/R}}
\left\{I^+_{E(\cW_R)}(p) + \varrho J_{V(\cW_R)}(p)\right\}
\end{equation*}
satisfies
\begin{equation}
\label{limvarformlb}
\liminf_{R\to\infty} \widehat{\chi}^{\,+}_R(\varrho) \geq \chi_{\cT}(\varrho).
\end{equation}
\end{itemize}
\end{proof}

%%%%%%%%% APPENDICES %%%%%%%%%%%%%%%%%%%%%%%%%%%%%%%%%

\appendix

%%%%%%%% APPENDIX A%%%%%%%%%%%%%%%%%%%%%%%%%%%%%%%%%%%

\section{Large deviation principle for the local times of Markov renewal processes}
\label{appA}

The following LDP, which was used in the proof of Lemma~\ref{lem:VF}, was derived in \cite[Proposition 1.2]{MZ2016}, and generalises the LDP for the empirical distribution of a Markov proceses on a finite state space derived in \cite{DV75}. See \cite[Chapter III]{dHLDP2000} for the definition of the LDP.
 
\begin{proposition}
\label{prop:LDPrenewal} 
Let $Y=(Y_t)_{t \geq 0}$ be the Markov renewal process on the finite graph $G=(V,E)$ with transition kernel $(\pi_{x,y})_{\{x,y\} \in E}$ and with sojourn times whose distributions $(\psi_x)_{x \in V}$ have support $(0,\infty)$. For $t > 0$, let $L_t^Y$ denote the empirical distribution of $Y$ at time $t$ (see \eqref{e:emp}). Then the family $(\mathbb{P}(L^Y_t \in\cdot))_{t>0}$ satisfies the LDP on $\mathcal{P}(V)$ with rate $t$ and with rate function $I^\dagger_E$ given by
\begin{equation*}
I^\dagger_E(p) = \inf_{\beta \in (0,\infty)} \inf_{q \in \cP(V)} \big[\widehat{K}(\beta q) + \widetilde{K}(p \mid \beta q)\big]
\end{equation*}
with 
\begin{eqnarray}
\label{rf1app}
\widehat{K}(\beta q) &=& \sup_{\widehat{q} \in \cP(V)} \sum_{x \in V} \beta q(x) \log\left(\tfrac{\widehat{q}(x)}{\sum_{y\in V} \pi_{x,y}\widehat{q}(y)}\right),\\    
\label{rf2app}
\widetilde{K}(p \mid \beta q ) &=& \sum_{x \in V} \beta q(x)\,(\cL\lambda_x)\left(\tfrac{p(x)}{\beta q(x)}\right),
\end{eqnarray}
where
\begin{equation*}
\begin{array}{llll}
(\cL\lambda_x)(\alpha) &=&  \sup_{\theta \in \mathbb{R}} [\alpha\theta - \lambda_x(\theta)], &\alpha \in [0,\infty),\\[0.2cm] 
\lambda_x(\theta) &=& \log \int_0^\infty \ee^{\theta \tau}\psi_x(\dd\tau), &\theta \in \mathbb{R}.
\end{array}
\end{equation*}
\end{proposition}

The rate function $I_E$ consist of two parts: $\widehat{K}$ in \eqref{rf1app} is the rate function of the LDP on $\cP(V)$ for the empirical distribution of the discrete-time Markov chain on $V$ with transition kernel $(\pi_{x,y})_{\{x,y\} \in E}$ (see \cite[Theorem IV.7]{dHLDP2000}), while $\widetilde{K}$ in \eqref{rf2app} is the rate function of the LDP on $\cP(0,\infty)$ for the empirical mean of the sojourn times, given the empirical distribution of the discrete-time Markov chain. Moreover, $\lambda_x$ is the cumulant generating function associated with $\psi_x$, and $\cL\lambda_x$ is the Legendre transform of $\lambda_x$, playing the role of the Cram\`er rate function for the empirical mean of the i.i.d.\ sojourn times at $x$. The parameter $\beta$ plays the role of the ratio between the continuous time scale and the discrete time scale.

%%%%%%%%%%% APPENDIX B %%%%%%%%%%%%%%%%%%%%%% 

\section{Sojourn times: cumulant generating functions and Legendre tranforms}
\label{appB}

In Appendix~\ref{app:general} we recall general properties of cumulant generating functions and Legendre transforms, in Appendices~\ref{app:exp} and \ref{app:nonexp} we identify both for the two sojourn time distributions arising in Lemma~\ref{lem:VF}, respectively.   

%%%

\subsection{General observations}
\label{app:general}

Let $\lambda$ be the cumulant generating function of a non-degenerate sojourn time distribution $\phi$, and $\cL\lambda$ be the Legendre transform of $\lambda$ (recall \eqref{cumulant}). Both $\lambda$ and $\cL\lambda$ are strictly convex, are analytic in the interior of their domain, and achieve a unique zero at $\theta = 0$, respectively, $\alpha=\alpha_c$ with $\alpha_c= \int_0^\infty \tau \phi(\dd\tau)$. Furthermore, $\lambda$ diverges at some $\theta_c \in (0,\infty]$ and has slope $\alpha_c$ at $\theta=0$. Moreover, if the slope of $\lambda$ diverges at $\theta_c$, then $\cL\lambda$ is finite on $(0,\infty)$. 

\begin{figure}[htbp]
\begin{center}
\setlength{\unitlength}{0.6cm}
\begin{picture}(12,8)(0,-.8)
%%%
\put(-.5,0){\line(6,0){9.5}}
\put(5,-2.5){\line(0,7){8.5}}
{\thicklines
\qbezier(0,-.5)(3,-.5)(5,0)
\qbezier(5,0)(7,0.6)(8,6)
}
\qbezier[60](4,-2)(6.4,1.2)(8.8,4.4)
\qbezier[20](6.4,0)(6.4,0.6)(6.4,1.2)
\qbezier[20](5,1.2)(6,1.2)(6.4,1.2)
\put(9.5,-0.2){$\theta$}
\put(4.3,6.5){$\lambda(\theta)$}
\put(9.2,4.75){$\alpha$}
\put(5.95,-.8){$\theta(\alpha)$}
\put(3.5,1.1){$\mu(\alpha)$}
\put(1.8,-1.1){$-(\cL\lambda)(\alpha)$}
\put(5,0){\circle*{.2}}
\put(6.4,1.2){\circle*{.2}}
\put(5,-.7){\circle*{.2}}
%%%
\end{picture}
%%%%%%%%%%%%%%%%%%%
\vspace{1cm}
\end{center}
\caption{Picture exhibiting the link between $\lambda(\theta)$, $(\cL\lambda)(\alpha)$, $\theta(\alpha)$, $\mu(\alpha)$. The dotted line is the tangent of $\lambda$ with slope $\alpha$, crossing the horizontal axis at $-(\cL\lambda)(\alpha)$, and touching $\lambda$ at the point $(\theta(\alpha),\mu(\alpha))$. All are analytic on the interior of their domain.}
\label{Fig:Legendre}
\end{figure}
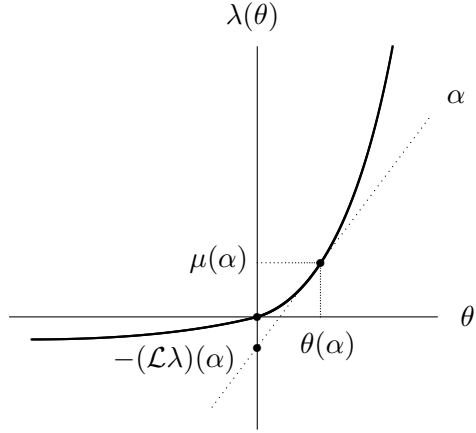

The supremum in the Legendre transform defining $(\cL\lambda)(\alpha)$ is uniquely taken at $\theta=\theta(\alpha)$ solving the equation 
\begin{equation*}
\lambda'(\theta(\alpha)) = \alpha.
\end{equation*}
The tangent of $\lambda$ with slope $\alpha$ at $\theta(\alpha)$ intersects the vertical axis at $(-\cL\lambda)(\alpha)$, i.e., putting 
\begin{equation}
\label{mulrel}
\mu(\alpha) = \lambda(\theta(\alpha))
\end{equation}
we have 
\begin{equation}
\label{mudef}
\mu(\alpha) = \alpha (\cL\lambda)'(\alpha)-(\cL\lambda)(\alpha).
\end{equation} 
(See Fig.~\ref{Fig:Legendre}.) Note that by differentiating \eqref{mudef} we get
\begin{equation*}
\mu'(\alpha) = \alpha(\cL\lambda)''(\alpha),
\end{equation*}   
which shows that $\alpha \mapsto \mu(\alpha)$ is strictly increasing and hence invertible, with inverse function $\mu^{-1}$.  
Note that by differentiating the relation $(\cL\lambda)(\alpha) = \alpha\theta(\alpha)-\lambda(\theta(\alpha))$ we get 
\begin{equation}
\label{cLrel1}
(\cL\lambda)'(\alpha) = \theta(\alpha).
\end{equation}
A further relation that is useful reads
\begin{equation}
\label{cLrel2}
(\cL\lambda)' \circ \mu^{-1} = \lambda^{-1},
\end{equation}
which follows because $\mu = \lambda \circ \theta$ by \eqref{mulrel} and $(\cL\lambda)' = \theta$ by \eqref{cLrel1}. 

%%%%%%%%%%%%%%%%%%%%%%%%%%%%%%%%%%%%%%%%%%%%%%%
%%%%%%%%%%%%%%%%%%%%%%%%%%%%%%%%%%%%%%%%%%%%%

%%%

\subsection{Exponential sojourn time}
\label{app:exp}

If $\phi=\mathrm{EXP}(d+1)$, then the cumulant generating function $\lambda(\theta) = \log \int_0^\infty \ee^{\theta\tau}\psi(\dd\tau)$ is given by
\begin{equation*}
\lambda(\theta) = 
\begin{cases}
\log\left(\tfrac{d+1}{d+1-\theta}\right), &\theta < d+1,\\[0.1cm]
\infty, &\theta \geq d+1.
\end{cases}
\end{equation*}
To find $(\cL\lambda)(\alpha)$, we compute 
\begin{equation*}
\frac{\partial}{\partial\theta}[\alpha\theta - \log(\tfrac{d+1}{d+1 - \theta})] = \alpha - \frac{1}{d+1-\theta},
\qquad \frac{\partial^2}{\partial\theta^2}[\alpha\theta - \log(\tfrac{d}{d+1-\theta})] = - \frac{1}{(d+1-\theta)^2} < 0.   
\end{equation*}
Hence the supremum in \eqref{Legendre} is uniquely taken at 
\begin{equation*}
\theta(\alpha) = d+1 - \tfrac{1}{\alpha}, \qquad \alpha > 0,
\end{equation*}
so that
\begin{equation}
\label{explaw}
(\cL\lambda)(\alpha) = \alpha (d+1) -1 - \log[\alpha (d+1)], \qquad \alpha>0.
\end{equation}
Thus, $\lambda$ and $\cL\lambda$ have the shape in Fig.~\ref{Fig:lambda}, with $\theta_c = d+1$ and $\alpha_c = \frac{1}{d+1}$, and with $\lim_{\theta \uparrow \theta_c} \lambda(\theta) = \infty$ and $\lim_{\theta \uparrow \theta_c} \lambda'(\theta) = \infty$.   

%%%%%%%%%%%%%%%%%%%%%%%%%%%%%%%%%%%%%%%%%%%%%%%
\begin{figure}[htbp]
\vspace{-1cm}
\begin{center}
\setlength{\unitlength}{0.45cm}
\begin{picture}(12,12)(7,-1)
%%%
\put(-.5,0){\line(6,0){9.5}}
\put(5,-1.5){\line(0,7){8.5}}
{\thicklines
\qbezier(0,-2)(3,-1.5)(5,0)
\qbezier(5,0)(6,0.8)(6.8,5)
}
\qbezier[40](2,-2.7)(5,0)(8,2.7)
\qbezier[40](7,0)(7,3)(7,6)
\put(4.4,0.4){$0$}
\put(9.5,-0.2){$\theta$}
\put(6.8,-.8){$\theta_c$}
\put(4.3,7.5){$\lambda(\theta)$}
\put(8,3){$\alpha_c$}
\put(5,0){\circle*{.2}}
%%%
\put(13.5,0){\line(6,0){9.5}}
\put(14,-1.5){\line(0,7){8.5}}
{\thicklines
\qbezier(14.3,6)(14.5,0)(17,0)
\qbezier(17,0)(20,0)(22,5)
}
\put(13.4,-.8){$0$}
\put(23.5,-0.2){$\alpha$}
\put(13.3,7.5){$(\cL\lambda)(\alpha)$}
\put(16.6,-.8){$\alpha_c$}
\put(17,0){\circle*{.2}}
%%%
\end{picture}
%%%%%%%%%%%%%%%%%%%
\vspace{0.5cm}
\end{center}
\caption{Picture of $\theta \mapsto \lambda(\theta)$ (left) and $\alpha \mapsto (\cL\lambda)(\alpha)$ (right) for $\phi=\mathrm{EXP}(d+1)$.}
\label{Fig:lambda}
\end{figure}
%%%%%%%%%%%%%%%%%%%%%%%%%%%%%%%%%%%%%%%%%%%%%%

\noindent
Note that $\mu$ has domain $(0,\infty)$ and range $\R$.

%%%

\subsection{Non-exponential sojourn time}
\label{app:nonexp}

For $\phi=\psi$ the computations are more involved. Let $\cT^*=(E,V)$ be the infinite rooted regular tree of degree $d+1$. Write $\cO$ for the root. Let $X = (X_n)_{n \in \N_0}$ be the discrete-time simple random walk on $\cT^*=(E,V)$ starting from $\cO$. Write $\tau_\cO$ to denote the time of the \emph{first return} of $X$ to $\cO$. Define $r = \mathbb{P}_\cO(\tau_\cO<\infty)$. It is easy to compute $r$ by projecting $X$ on $\N_0$: $r$ is the return probability to the origin of the random walk on $\N_0$ that jumps to the right with probability $p = \tfrac{d}{d+1}$ and to the left with probability $q = \tfrac{1}{d+1}$, which equals $\frac{p}{q}$ (see \cite[Section 8]{S1976}). Thus, $r= \frac{1}{d}$. 
  
For $y \in \cT^*$, define $h_y$ = $\mathbb{P}_y(\tau_\cO <\infty)$. Then $h_y$ can be explicitly calculated, namely,
\begin{equation*}
h_y = 
\begin{cases}
d^{-|y|},  &y\in \cT^*\setminus\{\cO\},\\
1, &y= \cO.
\end{cases}
\end{equation*}
Note that $h$ is a harmonic function on $\cT^* \setminus \cO$, i.e., $h_y = \sum_{z\in\cT^*} \widehat{\pi}_{y,z} h_z$, $y\in\cT^*\setminus \cO$. We can therefore consider the Doob-transform of $X$, which is the random walk with transition probabilities away from the root given by
\begin{equation*}
\check{\sigma}_{y,z} = 
\begin{cases} 
\frac{d}{d+1}, &z=y^\uparrow,\\[0.1cm]
\frac{1}{d}\frac{1}{d+1}, &z\neq y^\uparrow, \{y,z\} \in E,\\[0.1cm]
0, &\text{else},
\end{cases} \qquad \qquad y \in \cT^*\setminus\{\cO\},
\end{equation*}
and transition probabilities from the root are given by
\begin{equation*}
\check{\sigma}_{\cO,z} = 
\begin{cases}
\frac{1}{d}, &\{\cO,z\}\in E,\\
0, &\text{else}.
\end{cases}
\end{equation*}
Thus, the Doob-transform reverses the upward and the downward drift of $X$.

Recall from Lemma~\ref{lem:VF} that $\psi$ is the distribution of $\tau_\cO$ \emph{conditional} on $\{\tau_\cO<\infty\}$ and on $X$ leaving $\cO$ at time 0.
 
\begin{lemma}
Let $\lambda(\theta) = \log \int_0^\infty \ee^{\theta\tau}\psi(\dd\tau)$. Then
\begin{equation}
\label{lambdaid}
\ee^{\lambda(\theta)} 
= \begin{cases}
\frac{d+1-\theta}{2}\,\left[1- \sqrt{1- \tfrac{4d}{(d+1-\theta)^2}}\ \right], &\theta \in (-\infty,\theta_c],\\
\infty, &\text{else},
\end{cases}
\end{equation}
with $\theta_c = (\sqrt{d}-1)^2$. The range of $\exp \circ \lambda$ is $(0,\sqrt{d}\,]$, with the maximal value is uniquely taken at $\theta=\theta_c$.
\end{lemma}

\begin{proof}
To compute the moment-generating function of $\tau_\cO$, we consider the Doob-transform of $X$ and its projection onto $\mathbb{N}_0$. Let $p_{2k} = P(\tau_\cO = 2k)$. It is well-known that (see \cite[Section 8]{S1976})
\begin{equation}
\label{Green}
G^{p,q}(s) = \EE(s^{\tau_\cO} \mid \tau_\cO <\infty) = \sum_{k \in \N} s^{2k} p_{2k} = \frac{1}{2p}\left[1- \sqrt{1-4pqs^2}\right], \qquad |s| \leq 1.
\end{equation} 
Therefore we have 
\begin{equation}
\label{min}
\begin{aligned}
\ee^{\lambda(\theta)} = \EE(\ee^{\theta\tau_\cO}) 
&= \sum_{k \in \N} p_{2k}\, 
\left[\EE\left(\ee^{\theta\,\mathrm{EXP}(d+1)}\right)\right]^{2k-1}\\
&= \sum_{k \in \N} p_{2k} \left(\frac{d+1}{d+1 - \theta}\right)^{2k-1}
= \left(\frac{d+1 -\theta}{d+1}\right) G^{p,q}(s)
\end{aligned}
\end{equation} 
with
\begin{equation*}
p = \tfrac{1}{d+1}, \qquad q = \tfrac{d}{d+1}, \qquad s = \frac{d+1}{d+1-\theta}. 
\end{equation*}
Inserting \eqref{Green} into \eqref{min}, we get the formula for $\lambda(\theta)$. From the term in the square root we see that $\lambda(\theta)$ is finite if and only if $\theta \leq \theta_c = d+1-2\sqrt{d} = (\sqrt{d}-1)^2$. 
\end{proof}

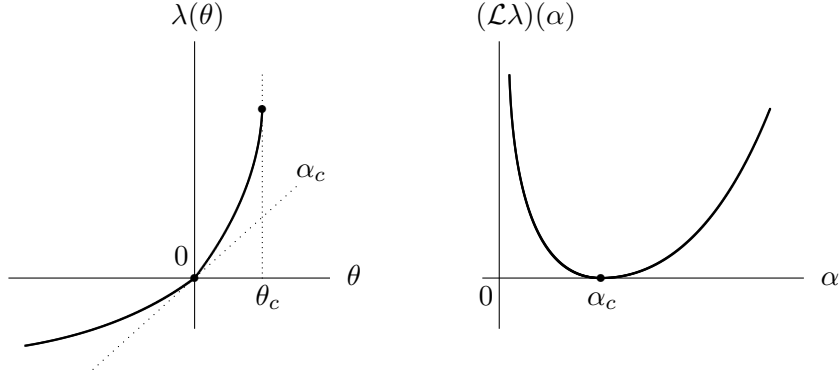
\begin{figure}[htbp]
\vspace{-1cm}
\begin{center}
\setlength{\unitlength}{0.45cm}
\begin{picture}(12,12)(7,-1)
%%%
\put(-.5,0){\line(6,0){9.5}}
\put(5,-1.5){\line(0,7){8.5}}
{\thicklines
\qbezier(0,-2)(3,-1.5)(5,0)
\qbezier(5,0)(6.9,2.5)(7,5)
}
\qbezier[40](2,-2.7)(5,0)(8,2.7)
\qbezier[40](7,0)(7,3)(7,6)
\put(4.4,0.4){$0$}
\put(9.5,-0.2){$\theta$}
\put(6.8,-.8){$\theta_c$}
\put(4.3,7.5){$\lambda(\theta)$}
\put(8,3){$\alpha_c$}
\put(5,0){\circle*{.2}}
%%%
\put(13.5,0){\line(6,0){9.5}}
\put(14,-1.5){\line(0,7){8.5}}
{\thicklines
\qbezier(14.3,6)(14.5,0)(17,0)
\qbezier(17,0)(20,0)(22,5)
}
\put(13.4,-.8){$0$}
\put(23.5,-0.2){$\alpha$}
\put(13.3,7.5){$(\cL\lambda)(\alpha)$}
\put(16.6,-.8){$\alpha_c$}
\put(17,0){\circle*{.2}}
\put(7,5){\circle*{.2}}
%%%
\end{picture}
%%%%%%%%%%%%%%%%%%%
\vspace{0.5cm}
\end{center}
\caption{Picture of $\theta \mapsto \lambda(\theta)$ (left) and $\alpha \mapsto (\cL\lambda)(\alpha)$ (right) for $\phi=\psi$.}
\label{Fig:lambdaalt}
\end{figure}

There is no easy closed form expression for $(\cL\lambda)(\alpha)$, but it is easily checked that $\lambda$ and $\cL\lambda$ have the shape in Fig.~\ref{Fig:lambdaalt}, with $\theta_c = (\sqrt{d}-1)^2$ and $\alpha_c = \int_0^\infty \tau \psi(\dd\tau)<\infty$, and with $\lambda(\theta_c) = \log \sqrt{d}<\infty$ and $\lambda'(\theta_c)=\infty$, i.e., there is a \emph{cusp} at the threshold $\theta_c$, implying that $\cL\lambda$ is finite on $(0,\infty)$. It follows from \eqref{cLrel1} that
\begin{equation}
\label{Legendrelim}
\lim_{\alpha \to \infty} \frac{1}{\alpha} (\cL\lambda)(\alpha) = \lim_{\alpha \to \infty} \theta(\alpha) = \theta_c.
\end{equation}

%%%%%%%%%%%%%%%%%%%%%%%%%%%%%%%%%%%%%%%%%%%%%%%
%%%%%%%%%%%%%%%%%%%%%%%%%%%%%%%%%%%%%%%%%%%%%%

\begin{lemma}
The function $\lambda^{-1} \circ \log = (\exp \circ \lambda)^{-1}$ is given by
\begin{equation}
\label{lambdainverse}
(\exp \circ \lambda)^{-1}(\beta) = d+1 - \beta -\frac{d}{\beta}, \qquad \beta \in (0,\sqrt{d}\,].
\end{equation}
The range of $(\exp \circ \lambda)^{-1}$ is $(-\infty,\theta_c]$, with the maximal value $\theta_c$ uniquely taken at $\beta = \sqrt{d}$.
\end{lemma}

\begin{proof}
We need to invert $\exp \circ \lambda$ in \eqref{lambdaid}. Abbreviate $\chi = \frac{d+1-\theta}{2}$. Then
\begin{equation*}
\beta = \chi \left[1-\sqrt{1-\frac{d}{\chi^2}}\,\right]  \quad \Longrightarrow \quad \chi = \frac{\beta^2+d}{2\beta} \quad
\Longrightarrow \quad \theta = d+1 - \frac{\beta^2 + d}{\beta}.   
\end{equation*}
\end{proof}

\noindent
Note that $(\sqrt{d},\infty)$ is not part of the domain of $(\exp \circ \lambda)^{-1}$, even though the right-hand side of \eqref{lambdainverse} still makes sense (as a second branch). Note that $\mu$ has domain $(0,\infty)$ and range $(-\infty,\sqrt{d}\,]$ (see Fig.~\ref{Fig:Legendre}).

%%%%%%%%% APPENDIX C %%%%%%%%%%%%%%%%%%%%%%%%%%%%%%

\section{Analysis of the variational problem on the infinite regular tree}
\label{appC}

In this appendix we prove Theorem~\ref{t:chivar}. Appendix~\ref{properties} formulates two theorems that imply Theorem~\ref{t:chivar}, Appendix~\ref{profile} provides the proof of these theorems. Recall the definition of  $\cP(V)$, $I_E(p)$ and $J_V(p)$ from \eqref{e:defIJ}. Set
\begin{equation}
\label{chivar}
\chi_\cT(\varrho) = \inf_{p \in \cP_\cO(V)} [I_E(p) + \varrho J_V(p)], \qquad \varrho \in (0,\infty),
\end{equation}
where $\cP_\cO(V) = \{p \in \cP(V)\colon\, \mathrm{argmax}\,p = \cO\}$. Since $\cP(V)$, $I_E$ and $J_V$ are invariant under translations, the centering at $\cO$ is harmless. 

%%%

\subsection{Two properties}
\label{properties}

\begin{theorem}
\label{thm:1}
For every $\varrho \in (0,\infty)$ the infimum in \eqref{chivar} is attained, and every minimiser $\Bar{p}$ is strictly positive, non-increasing in the distance to the root, and such that 
\begin{equation*}
\sum_{N\in\N_0} \partial S_R \log (R+1) \leq \frac{d+1}{\varrho},
\qquad \partial S_R = \sum_{\partial B_R(\cO)} \Bar{p}(x), 
\end{equation*}
where $B_R(\cO)$ is the ball of radius $R$ around $\cO$.
\end{theorem}

\begin{theorem}
\label{thm:2}
The function $\varrho \mapsto \chi_\cT(\varrho)$ is strictly increasing and globally Lipschitz continuous on $(0,\infty)$, with $\lim_{\varrho \downarrow 0} \chi_\cT(\varrho) = d-1$ and $\lim_{\varrho \to \infty} \chi_\cT(\varrho) = d+1$. 
\end{theorem}

Theorems~\ref{thm:1}--\ref{thm:2} settle Theorem~\ref{t:chivar}. Their proof uses the following two  lemmas.

\begin{lemma}
\label{lem:Jupper}
For every $\varrho \in (0,\infty)$, the infimum in \eqref{chivar} may be restricted to $p \in \cP_\cO(V)$ such that $J_V(p) \leq \tfrac{d+1}{\varrho}$.    
\end{lemma}

\begin{proof}
Let $\delta_\cO \in \cP_\cO(V)$ denote the point measure at $\cO$. Then, for all $\varrho \in (0,\infty)$, 
\begin{equation*}
\chi_\cT(\varrho) \leq I_E(\delta_\cO) + \varrho J_V(\delta_\cO) =  (d+1) + \varrho \times 0 = d+1.   
\end{equation*}
Since $I_V \geq 0$, we may restrict the infimum in \eqref{chivar} to $p$ with $J_V(p) \leq \frac{d+1}{\varrho}$.
\end{proof}

\begin{lemma}
\label{lem:Jlower}
For every $\varrho \in (0,\infty)$, there exists a $c(\varrho) >0$ such that the infimum in \eqref{chivar} may be restricted to $p\in\mathcal{P}_\cO(V)$ such that $J_V(p) \geq c(\varrho)$.  
\end{lemma}

\begin{proof}
Since $J_V(p) = 0$ if and only if $p = \delta_\cO$ is a point measure, it suffices to show that $\delta_\cO$ is not a minimiser of $\chi_\cT(\varrho)$. To that end, for $y \in V$ compute
\begin{equation}
\label{Lagrange}
\frac{\partial}{\partial p(y)}[I_E(p) + \varrho J_V(p)] = 1 - \sum_{z\sim y}\sqrt{\frac{p(z)}{p(y)}} - \varrho\log p(y) -\varrho.     
\end{equation}
Because $p(\cO)>0$, it follows that the right-hand side tends to $-\infty$ as $p(y) \downarrow 0$ for every $y \sim \cO$. Hence, no $p \in \cP_\cO(V)$ with $p(y) = 0$ for some $y \sim \cO$ can be a minimiser of \eqref{chivar}, or be the weak limit point of a minimising sequence. In particular, $\delta_\cO$ cannot.
\end{proof}
 
%%%

\subsection{Proof of the two properties}
\label{profile}

\begin{proof}[Proof of Theorem~\ref{thm:1}]
First observe that $\cP(V)$ and $J_V$ are invariant under permutations, i.e., for any $p \in \cP(V)$ and any relabelling $\pi$ of the vertices in $V$, we have $\pi p \in \cP(V)$ and $J_V(\pi p)=J_V(p)$. The same does not hold for $I_E$, but we can apply permutations such that $I_E(\pi p) \leq I_E(p)$. 

\medskip\noindent
{\bf 1.}
Pick any $p \in \cP(V)$. Pick any backbone $\bb = \{x_0, x_1,\cdots\}$ that runs from $x_0 = \cO$ to infinity. Consider a permutation $\pi$ that \emph{reorders} the vertices in $\bb$ such that $\{(\pi p)(x)\}_{x \in \bb}$ becomes \emph{non-increasing}. Together with the reordering, transport all the trees that hang off $\bb$ as well. Since $\pi p$ is non-increasing along $\bb$, while all the edges that do not lie on $\bb$ have the same neighbouring values in $p$ and in $\pi p$, we have \begin{equation}
\label{order}
I_E(\pi p) \leq I_E(p).
\end{equation}
Indeed, 
\begin{equation}
\label{diff}
\tfrac12\, [I_E(p) - I_E(\pi p)] =  \sum_{k \in \N_0} \sqrt{(\pi p)(x_k) (\pi p)(x_{k+1})} 
- \sum_{k \in \N_0} \sqrt{p(x_k)p(x_{k+1})},
\end{equation}
where we use that $p(x_0) = (\pi p)(x_0)$ (because $p(x_0) \geq p(x_k)$ for all $k\in\N$) and $\sum_{k\in\N} p(x_k) = \sum_{k\in\N} (\pi p)(x_k)$. The right-hand side of \eqref{diff} is $\geq 0$ by the rearrangement inequality for sums of products of two sequences \cite[Section 10.2, Theorem 368]{HLP1952}. In fact, strict inequality in \eqref{diff} holds unless $p$ is constant along $\bb$. But this is impossible possible because it would imply that $p(\cO) = 0$ and hence $p(x) = 0$ for all $x \in V$. Thus, $p$ and $\bb$ being arbitrary, it follows from \eqref{order} that any minimiser or minimising sequence must be non-increasing in the distance to $\cO$.  Indeed, if it were not, then there would be a $\bb$ along which the reordering would lead to a lower value of $I_E+\varrho J_V$. Hence we may replace \eqref{chivar} by
\begin{equation}
\label{chivaralt}
\chi_\cT(\varrho) = \inf_{p \in \cP_\cO^\downarrow(V)} [I_E(p) + \varrho J_V(p)], \qquad \varrho \in (0,\infty),
\end{equation}
with $\cP_\cO^\downarrow(V)$ defined in \eqref{pdownset}.

\medskip\noindent
{\bf 2.}
Let $p \in \cP_\cO^\downarrow(V)$. Estimate
\begin{equation*}
J_V(p) = \sum_{R \in \N_0}  \sum_{x \in \partial B_R(\cO)} [-p(x)\log p(x)]
\geq \sum_{R \in \N_0}\sum_{x \in \partial B_R(\cO)} \Big[-p(x)\log\big(\tfrac{1}{R+1}\big)\Big],
\end{equation*}
where we use that $p(x) \leq \tfrac{1}{R+1}$ for all $x \in \partial B_R(\cO)$. Hence
\begin{equation*}
J_V(p) \geq \sum_{R \in \N_0} \partial S_R \log(R+1)
\end{equation*}
with $\partial S_R = \sum_{x \in \partial B_R(\cO)} p(x)$. By Lemma~\ref{lem:Jupper}, $J_V(p) \leq \frac{d+1}{\varrho}$, and so 
\begin{equation}
\label{eq:tightness}
\sum_{R \in \N_0} \partial S_R \log(R+1) \leq \frac{d+1}{\varrho}.
\end{equation}
The computation in \eqref{Lagrange} shows that any $p$ for which there exist $z \sim y$ with $p(z)>0$ and $p(y)=0$ cannot be minimiser nor a weak limit point of a minimising sequence. Hence all minimisers or weak limit points of minimising sequences are strictly positive everywhere.    

\medskip\noindent
{\bf 3.} 
Take any minimising sequence $(p_n)_{n\in\N}$ of \eqref{chivaralt}. By \eqref{eq:tightness}, $\lim_{R\to\infty} \sum_{x \notin B_R(\cO)} p_n(x) = 0$ uniformly in $n\in\N$, and so $(p_n)_{n\in\N}$ is tight. By Prokhorov's theorem, tightness is equivalent to $(p_n)_{n\in\N}$ being relatively compact, i.e., there is a subsequence $(p_{n_k})_{k\in\N}$ that converges weakly to a limit $\Bar{p} \in \cP_\cO^\downarrow(V)$. By Fatou's lemma, we have $\liminf_{k\to\infty} I_E(p_{n_k}) \geq I_E(\Bar{p})$ and $\liminf_{k\to\infty} J_V(p_{n_k}) \geq J_V(\Bar{p})$. Hence
\begin{equation*}
\chi_\cT(\varrho) = \lim_{k \to \infty} [I_E(p_{n_k}) + \varrho J_V(p_{n_k})] \geq I_E(\Bar{p}) + \varrho J_V(\Bar{p}). 
\end{equation*}
Hence $\Bar{p}$ is a minimiser of \eqref{chivaralt}.
\end{proof} 

\begin{proof}[Proof of Theorem~\ref{thm:2}]
The proof uses approximation arguments.

\medskip\noindent
{\bf 1.}
We first show that $\varrho \mapsto \chi_\cT(\varrho)$ is strictly increasing and globally Lipschitz. Pick $\varrho_1 < \varrho_2$. Let $\bar{p}_{\varrho_1}$ be any minimiser of \eqref{chivar} at $\varrho_1$, i.e.,
\begin{equation*}
\chi_\cT(\varrho_1) = I_E(\bar{p}_{\varrho_1}) + \varrho_1 J_V(\bar{p}_{\varrho_1}).
\end{equation*}
Estimate
\begin{equation*}
\begin{aligned} 
&[I_E(\bar{p}_{\varrho_1}) + \varrho_1 J_V(\bar{p}_{\varrho_1})] 
= [I_E(\bar{p}_{\varrho_1}) + \varrho_2 J_V(\bar{p}_{\varrho_1})] - (\varrho_2 - \varrho_1)J_V(\bar{p}_{\varrho_1})\\
&\geq \chi_\cT(\varrho_2) - (\varrho_2 - \varrho_1) J_V(\bar{p}_{\varrho_1}) 
\geq \chi(\varrho_2) - (\varrho_2 - \varrho_1) \tfrac{d+1}{\varrho_1},
\end{aligned}
\end{equation*}
where we use Lemma~\ref{lem:Jupper}. Therefore 
\begin{equation*}
\chi_\cT(\varrho_2) - \chi_\cT(\varrho_1) \leq (\varrho_2-\varrho_1) \tfrac{d+1}{\varrho_1}.
\end{equation*}
Similarly, let $\bar{p}_{\varrho_2}$ be any minimiser of \eqref{chivar} at $\varrho_2$, i.e.,
\begin{equation*}
\chi_\cT(\varrho_2) = I_E(\bar{p}_{\varrho_2}) + \varrho_2 J_V(\bar{p}_{\varrho_2}).
\end{equation*}
Estimate
\begin{equation*}
\begin{aligned}
&[I_E(\bar{p}_{\varrho_2}) + \varrho_2 J_V(\bar{p}_{\varrho_2})] 
= [I_E(\bar{p}_{\varrho_2}) + \varrho_1 J_V(\bar{p}_{\varrho_2})] + (\varrho_2 - \varrho_1) J_V(\bar{p}_{\varrho_2})\\
&\geq \chi_\cT(\varrho_1) + (\varrho_2 - \varrho_1) J_V(\bar{p}_{\varrho_2}) 
\geq \chi_\cT(\varrho_1) + (\varrho_2 - \varrho_1) c(\varrho_2),
\end{aligned}
\end{equation*}
where we use Lemma~\ref{lem:Jlower}. Therefore
\begin{equation*}
\chi_\cT(\varrho_2) - \chi_\cT(\varrho_1) \geq c(\varrho_2)(\varrho_2 - \varrho_1).     
\end{equation*}

\medskip\noindent
{\bf 2.}
Because $\chi_\cT(\varrho) \leq d+1$ for all $\varrho \in (0,\infty)$, it follows that $\lim_{\varrho \to \infty} \chi_\cT(\varrho) \leq d+1$. To obtain the reverse inequality, let $\Bar{p}_\varrho$ be any minimiser of \eqref{chivaralt} at $\varrho$. By Lemma~\ref{lem:Jupper}, we may assume that $J_V(\Bar{p}_\varrho) \leq \frac{d+1}{\varrho}$. Hence $\lim_{\varrho \to \infty} J_V(\Bar{p}_\varrho) = 0$, and consequently $\lim_{\varrho \to \infty} \Bar{p}_\varrho= \delta_\cO$ weakly. Therefore, by Fatou's lemma, $\lim_{\varrho \to \infty} \chi_\cT(\varrho) = \lim_{\varrho \to \infty} [I_E(\Bar{p}) + \varrho J_V(\Bar{p})] \geq \liminf_{\varrho \to \infty} I_E(\Bar{p}_\varrho) \geq I_E(\delta_\cO) = d+1$. 

\medskip\noindent
{\bf 3.} 
To prove that $\lim_{\varrho \downarrow 0}\chi_\cT(\varrho) \leq d-1$, estimate
\begin{equation*}
\chi_\cT(\varrho) \leq \inf_{ {p \in \cP_\cO^\downarrow(V)} \atop {\supp(p) \subseteq B_R(\cO)} } [I_E(p)+\varrho J_V(p)],
\qquad R \in \N_0.
\end{equation*}
Because 
\begin{equation*}
\sup_{ {p \in \cP_\cO^\downarrow(V)} \atop {\supp(p) \subseteq B_R(\cO)} } J_V(p) = J_V(p_R) = \log |B_R(\cO)|,
\qquad R \in \N_0,
\end{equation*}
with
\begin{equation*}
p_R(x) = 
\begin{cases}
|B_R(\cO)|^{-1}, &x \in B_R(\cO),\\[0.1cm]
0, &\text{else},
\end{cases}  
\end{equation*}
it follows that
\begin{equation*}
\lim_{\varrho \downarrow 0} \chi_\cT(\varrho) 
\leq \inf_{ {p \in \cP_\cO^\downarrow(V)} \atop {\supp(p) \subseteq B_R(\cO)} } I_E(p)
\leq I_E(p_R), \qquad R \in \N_0.
\end{equation*}
Compute (recall \eqref{e:defIJ}) ,
\begin{equation*}
I_E(p_R) = \frac{|\partial B_{R+1}(\cO)|}{|B_R(\cO)|}, \qquad R \in \N_0. 
\end{equation*}
Inserting the relations
\begin{equation*}
\begin{aligned}
&|\partial B_{R}(\cO)| = \left\{\begin{array}{ll}
1, &R=0,\\
(d+1)d^{R-1}, &R \in \N,
\end{array}
\right.
\\
&|B_R(\cO)| = \sum_{R'=0}^R |\partial B_{R'}(\cO)| = 1 + \frac{d+1}{d-1}(d^R-1), 
\quad R \in \N_0,
\end{aligned}
\end{equation*}
we get
\begin{equation*}
\begin{aligned}
I_E(p_R) = (d-1)\,\frac{(d+1)d^R}{(d+1)d^R-2}. 
\end{aligned}
\end{equation*}
Hence $\lim_{R\to\infty} I_E(p_R) = d-1$, and so $\lim_{\varrho \downarrow 0} \chi_\cT(\varrho) \leq d-1$. 

\medskip\noindent
{\bf 4.}
To prove that $\lim_{\varrho \downarrow 0} \chi_\cT(\varrho) \geq d-1$, note that because $J_V \geq 0$ we can estimate 
\begin{equation*}
\lim_{\varrho \downarrow 0} \chi_\cT(\varrho) \geq \inf_{p \in \cP_\cO^\downarrow(V)} I_E(p).
\end{equation*}
It therefore suffices to show that 
\begin{equation*}
\inf_{p \in \cP_\cO^\downarrow(V)} I_E(p) \geq d-1, 
\end{equation*}
i.e., $(p_R)_{R \in \N_0}$ is a minimising sequence of the infimum in the left-hand side. The proof goes as follows. Write (recall \eqref{e:defIJ})
\begin{equation*}
\begin{aligned}
I_E(p) &= \tfrac12 \sum_{ {x,y \in V} \atop {x \sim y} } \left(\sqrt{p(x)} - \sqrt{p(y)}\,\right)^2\\
&= \tfrac12 \sum_{ {x,y \in V} \atop {x \sim y} } \left[p(x) + p(y) - 2 \sqrt{p(x)p(y)}\,\right]
= (d+1) - \sum_{ {x,y \in V} \atop {x \sim y} } \sqrt{p(x)p(y)}.
\end{aligned}
\end{equation*}
Since $\cT$ is a tree, each edge can be labelled by the end-vertex that is farthest from $\cO$. Hence the sum in the right-hand side can be written as
\begin{equation*}
\sum_{x \in V \setminus \cO } 2\sqrt{p(x)p(x^\downarrow)},
\end{equation*}
where $x^\downarrow$ is the unique neighbour of $x$ that is closer to $\cO$ than $x$. Since $2\sqrt{p(x)p(x^\downarrow)} \leq p(x) + p(x^\downarrow)$, it follows that 
\begin{equation*}
\sum_{x \in V \setminus \cO} 2\sqrt{p(x)p(x^\downarrow)} 
\leq \sum_{x \in V \setminus \cO} p(x) + \sum_{x \in V \setminus \cO} p(x^\downarrow) 
= [1-p(\cO)] + 1.
\end{equation*}
Therefore
\begin{equation*}
I_E(p) \geq d - 1 + p(\cO),
\end{equation*}
which settles the claim.
\end{proof}

%%%%%%%%%%%%%%%%%%%%%%%%%%%%%%%%%%%%%%

\section{Large deviation estimate for the local time away from the backbone}
\label{appD}

In this appendix we derive a large deviation principle for the \emph{total local times at successive depths} of the random walk on $\cT^\Z$ (see Fig.~\ref{fig:redrawing}). This large deviation principle is not actually needed, but serves as a warm up for the more elaborate computations in Appendix~\ref{appE}. 
  
For $k\in\N_0$, let $V_k$ be the set of vertices in $\cT^\Z$ that are at distance $k$ from the backbone (see Fig.~\ref{fig:redrawing}). For $R \in \N$, define 
\begin{equation*}
\begin{array}{llll}
\ell^R_t(k) &=& \sum_{x \in V_k} \ell^\Z_t(x),  &k = 0,1,\ldots,R,\\[0.1cm]
\ell_t^R &=& \sum_{k > R} \sum_{x\in V_k} \ell^\Z_t(x), &k= R+1,
\end{array}
\end{equation*}
and
\begin{equation*}
L_t^R = \frac{1}{t}\,\Big((\ell_t(k))_{k=0}^R, \ell^R_t\Big).
\end{equation*}
Abbreviate $V^*_R = \{0,1,\ldots,R,R+1\}$,

\begin{lemma}
For every $R \in \N$, $(L_t^R)_{t \geq 0}$ satisfies the large deviation principle on $\cP(V^*_R)$ with rate $t$ and with rate function $I^\dagger_R$ given by
\begin{equation}
\label{IRid}
\begin{aligned}
I^\dagger_R(p) &= \big[\sqrt{(d-1)p(0)}-\sqrt{dp(1)}\,\big]^2 + \sum_{k=1}^{R-1} \big[\sqrt{p(k)}-\sqrt{dp(k+1)}\,\big]^2\\
&\qquad + \big[\sqrt{p(R)+p(R+1)} - \sqrt{dp(R+1)}\,\big]^2.
\end{aligned}
\end{equation}
\end{lemma}

\begin{proof}
By monitoring the random walk on the tree in Fig.~\ref{fig:redrawing} and projecting its depth on the vertices $0,1,\ldots,R$, respectively, $R+1$, we can apply the LDP in Proposition~\ref{prop:LDPrenewal} (see Fig.~\ref{fig:depthprojection}).

%%%%%%%%%%%%%%%%%%%%%%%%%%%%%%%%%%%%%%%%%%%%%
\begin{figure}[htbp]
\begin{center}
\setlength{\unitlength}{1cm}
\begin{picture}(8,1)(1,0)
%%%
\put(1,0){\line(1,0){7}}
%%%
\put(.9,-.7){$0$}
\put(1.9,-.7){$1$}
\put(6.85,-.7){$R$}
\put(7.65,-.7){$R+1$}
\put(1,0){\circle*{.2}}
\put(2,0){\circle*{.2}}
\put(7,0){\circle*{.2}}
\put(7.95,-.1){$\Box$}
\put(4,-.7){$\cdots$}
%%%
\end{picture}
\vspace{1cm}
\end{center}
\caption{\small Depths $k=0,1,\ldots,R$ and $k>R$.}
\label{fig:depthprojection}
\end{figure}
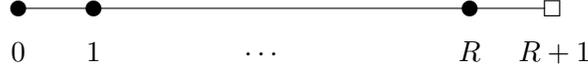
%%%%%%%%%%%%%%%%%%%%%%%%%%%%%%%%%%%%%%%%%%%%%%%%%

\medskip\noindent
{\bf 1.}
The sojourn times have distribution $\mathrm{EXP}(d+1)$ at vertices $k=0,1,\ldots,R$ and distribution $\psi$ at vertex $k=R+1$.  The transition probabilities are 
\begin{equation*}
\begin{array}{llll}
&\pi_{0,0} = \tfrac{2}{d+1}, &\pi_{0,1} = \tfrac{d-1}{d+1}, &\\[0.2cm]
&\pi_{k,k+1} = \tfrac{1}{d+1}, &\pi_{k,k-1} = \tfrac{d}{d+1}, &\quad k = 1,\ldots,R,\\[0.2cm]
&\pi_{R+1,R} = 1. &&
\end{array}
\end{equation*}
Proposition~\ref{prop:LDPrenewal} therefore yields that $(L_t^R)_{t \geq 0}$ satisfies the LDP on on $\cP(V^*_R)$ with rate $t$ and with rate function $I^\dagger_R$ given by
\begin{equation}
\label{IEdef}
I^\dagger_R(p) = (d+1) \sum_{k=0}^R p(k) + \inf_{v\colon V^*_R \to (0,\infty)} \sup_{u\colon V^*_R \to (0,\infty)} L(u,v)  
\end{equation}
with 
\begin{equation}
\label{Luv}
L(u,v) = - A -  B - C,
\end{equation} 
where 
\begin{equation*}
\begin{aligned}
A &= \sum_{k=1}^R v(x) \left\{1+\log\left(\frac{du(k-1)+u(k+1)}{u(k)}\,\frac{p(k)}{v(k)}\right)\right\},\\
B &= v(0) \left\{1+\log\left(\frac{2u(0)+(d-1)u(1)}{u(0)}\,\frac{p(0)}{v(0)}\right)\right\},\\
C &= v(R+1) \left\{\log\left(\frac{u(R)}{u(R+1)}\right)-(\cL\lambda)\left(\frac{p(R+1)}{v(R+1)}\right)\right\}.
\end{aligned}
\end{equation*}
Here we use \eqref{explaw} to compute $A$ and $B$, and for $C$ we recall that $\cL\lambda$ is the Legendre transform of the cumulant generation function $\lambda$ of $\psi$ computed in Lemma~\ref{lambdaid}.

\medskip\noindent
{\bf 2.}
We compute the infimum of $L(u,v)$ over $v$ for fixed $u$.

\medskip\noindent
$\bullet$ For $k=1,\ldots,R$,
\begin{equation*}
\begin{aligned}
&\frac{\partial A}{\partial v(k)} = \log\left(\frac{du(k-1)+u(k+1)}{u(k)}\,\frac{p(k)}{v(k)}\right),\\
&\Longrightarrow \bar{v}_u(k) = p(k)\,\frac{du(k-1)+u(k+1)}{u(k)}.
\end{aligned}
\end{equation*}
The second derivative is $1/v(k)>0$. 

\medskip\noindent
$\bullet$ For $k=0$,
\begin{equation*}
\begin{aligned}
&\frac{\partial B}{\partial v(0)} = \log\left(\frac{2u(0)+(d-1)u(1)}{u(0)}\,\frac{p(0)}{v(0)}\right),\\
&\Longrightarrow \bar{v}_u(0) = p(0)\,\frac{2u(0)+(d-1)u(1)}{u(0)}.
\end{aligned}
\end{equation*}
The second derivative is $1/v(0)>0$. 

\medskip\noindent
$\bullet$ For $k=R+1$, the computation is more delicate. Define (recall \eqref{mudef} in Appendix~\ref{appB})
\begin{equation*}
\mu(\alpha) = \alpha (\cL\lambda)^{'}(\alpha) - (\cL\lambda)(\alpha).
\end{equation*}
The function $\mu$ has range $(-\infty,\log \sqrt{d}\,]$, with the maximal value uniquely taken at $\alpha=\infty$. Therefore there are two cases.

\medskip\noindent
$\blacktriangleright$ $u(R+1)/u(R) \leq \sqrt{d}$. Compute
\begin{equation*}
\begin{aligned}
&\frac{\partial C}{\partial v(R+1)} = \mu\left(\frac{p(R+1)}{v(R+1)}\right) - \log\left(\frac{u(R+1)}{u(R)}\right),\\
&\Longrightarrow \bar{v}(R+1) = \frac{p(R+1)}{\alpha_u(R+1)}
\end{aligned} 
\end{equation*}
with $\alpha_u(R+1)$ solving the equation
\begin{equation*}
\log \left(\frac{u(R+1)}{u(R)}\right) = \mu\big(\alpha_u(R+1)\big).
\end{equation*}
Since $\mu'(\alpha) = \alpha(\cL\lambda)''(\alpha)$ and $\cL\lambda$ is strictly convex (see Fig.~\ref{Fig:lambdaalt} in Appendix \ref{appB}), $\mu$ is strictly increasing and therefore invertible. Consequently,
\begin{equation}
\label{alphaid}
\alpha_u(R+1) = \mu^{-1}\left(\log\left(\frac{u(R+1)}{u(R)}\right)\right).
\end{equation} 
Putting \eqref{Luv}--\eqref{alphaid} together, we get
\begin{equation}
\label{Lu}
L(u) = \inf_{v\colon V^*_R \to (0,\infty)} L(u,v)\\ 
= - \sum_{k=1}^R A_u(k) - B_u + C_u
\end{equation}
with
\begin{equation*}
\begin{aligned}
A_u(k) &= \frac{du(k-1)+u(k+1)}{u(k)}\,p(k), \qquad k = 1,\ldots,R,\\
B_u &= \frac{2u(0)+(d-1)u(1)}{u(0)}\,p(0),
\end{aligned}
\end{equation*}
and 
\begin{equation*}
\begin{aligned}
C_u &= \frac{p(R+1)}{\alpha_u(R+1)}\left[(\cL\lambda)\big(\alpha_u(R+1)\big) - \log\left(\frac{u(R+1)}{u(R)}\right)\right]\\
&= \frac{p(R+1)}{\alpha_u(R+1)} \big[(\cL\lambda)\big(\alpha_u(R+1)\big) - \mu\big(\alpha_u(R+1)\big)\big]\\ 
&= p(R+1)\,(\cL\lambda)^{'}\big(\alpha_u(R+1)\big)\\
&= p(R+1)\,((\cL\lambda)^{'} \circ \mu^{-1})\left(\log \left(\frac{u(R+1)}{u(R)}\right)\right).
\end{aligned}
\end{equation*}
In \eqref{cLrel2} in Appendix~\ref{appB} we showed that $(\cL\lambda)' \circ \mu^{-1} = \lambda^{-1}$. Moreover, in \eqref{lambdainverse} in Appendix~\ref{appB} we showed that $(\lambda^{-1} \circ \log) = S$ with
\begin{equation}
\label{Sdef}
S(\beta) = d+1 - \beta - \frac{d}{\beta}, \qquad \beta \in (0,\sqrt{d}\,].
\end{equation} 
Since $S$ has domain $(0,\sqrt{d}\,]$, $C_u(R+1)$ is only defined when $u(R+1)/u(R) \leq \sqrt{d}$, in which case
\begin{equation}
\label{Cu}
C_u = p(R+1)\,S\left(\frac{u(R+1)}{u(R)}\right).
\end{equation}

\medskip\noindent
$\blacktriangleright$ $u(R+1)/u(R) \leq \sqrt{d}$. In this case $\frac{\partial C}{\partial v(R+1)}>0$, the infimum is taken at $\bar{v}(R+1)=0$, and hence (recall \eqref{Legendrelim})
\begin{equation}
\label{Cualt}
C_u = p(R+1)\,(\sqrt{d}-1)^2 = p(R+1)\, S(\sqrt{d}). 
\end{equation}
Note that the right-hand side does not depend on $u$. The expressions in \eqref{Cu}--\eqref{Cualt} can be summarised as
\begin{equation*}
C_u = p(R+1)\,S\left(\sqrt{d} \wedge \frac{u(R+1)}{u(R)}\right).
\end{equation*} 

\medskip\noindent
{\bf 3.}
Next we compute the supremum over $u$ of 
\begin{equation}
\label{Lubarv}
L(u) = L(u,\bar{v}_u) = - A_u - B_u + C_u.
\end{equation} 
with $A_u = \sum_{k=1}^R A_u(k)$. We only write down the derivatives that are non-zero. 

\noindent
$\bullet$ For $k=2,\ldots,R-1$,
\begin{equation*}
- \frac{\partial A_u}{\partial u(k)} = - p(k+1)\,\frac{d}{u(k+1)} - p(k-1)\,\frac{1}{u(k-1)} + p(k)\,\frac{du(k-1)+u(k+1)}{u(k)^2}.
\end{equation*} 
$\bullet$ For $k=1$,
\begin{equation*}
\begin{aligned}
&- \frac{\partial A_u}{\partial u(1)} = - p(2)\,\frac{d}{u(2)} + p(1)\,\frac{du(0)+u(2)}{u(1)^2},\\
&- \frac{\partial B_u}{\partial u(1)} = - p(0)\,\frac{d-1}{u(0)}.
\end{aligned}
\end{equation*}
$\bullet$ For $k=R$,
\begin{equation*}
\begin{aligned}
&- \frac{\partial A_u}{\partial u(R)} = - p(R-1)\,\frac{1}{u(R-1)} + p(R)\,\frac{du(R-1)+u(R+1)}{u(R)^2},\\
&\frac{\partial C_u}{\partial u(R)} = p(R+1)\,\left[\frac{u(R+1)}{u(R)^2} - \frac{d}{u(R+1)}\right]\,
1_{\big\{\frac{u(R+1)}{u(R)} \leq \sqrt{d}\big\}}.
\end{aligned}
\end{equation*}
$\bullet$ For $k=0$,
\begin{equation*}
\begin{aligned}
&-\frac{\partial A_u}{\partial u(0)} = - p(1)\,\frac{d}{u(1)},\\
&-\frac{\partial B_u}{\partial u(0)} = p(0)\,\frac{(d-1)u(1)}{u(0)^2}.   
\end{aligned}
\end{equation*}
$\bullet$ For $k=R+1$,
\begin{equation*}
\begin{aligned}
&-\frac{\partial A_u}{\partial u(R+1)} = - p(R)\,\frac{1}{u(R)},\\
&\frac{\partial C_u}{\partial u(R+1)} = p(R+1)\,\left[-\frac{1}{u(R)} + \frac{du(R)}{u(R+1)^2}\right]\,
1_{\big\{\frac{u(R+1)}{u(R)} \leq \sqrt{d}\big\}}.
\end{aligned} 
\end{equation*}
All the first derivatives of $A_u+B_u+C_u$ are zero when we choose
\begin{equation}
\label{ubar}
\begin{aligned}
&\bar{u}(0) = \sqrt{(d-1)p(0)}, \qquad \bar{u}(k) = \sqrt{d^kp(k)}, \quad k = 1,\ldots,R,\\
&\bar{u}(R+1) = \sqrt{d^{R+1}\,\frac{p(R)p(R+1)}{p(R)+p(R+1)}}.
\end{aligned}
\end{equation} 
All the second derivatives are strictly negative, and so $\bar{u}$ is the unique maximiser. 

\medskip\noindent
{\bf 4.}
Inserting \eqref{ubar} into \eqref{Lu}, we get
\begin{equation*}
\begin{aligned}
&L(\bar{u}) = L(\bar{u},\bar{v}_{\bar{u}}) = - \sum_{k=2}^{R-1} A_{\bar{u}}(k) 
- \big[A_{\bar{u}}(1) + B_{\bar{u}}\big] - A_{\bar{u}}(R) + C_{\bar{u}}\\
&= -\sum_{k=2}^{R-1} \sqrt{dp(k)}\,\big[\sqrt{p(k-1)} + \sqrt{p(k+1)}\,\big]\\
&\qquad - \big[2\sqrt{d(d-1)p(0)p(1)} + 2p(0) + \sqrt{dp(1)p(2)}\,\big]\\
&\qquad  - \left[\sqrt{dp(R-1)p(R)} + \sqrt{\frac{p(R)}{p(R)+p(R+1)}}\,\sqrt{dp(R)p(R+1)}\,\right]\\
&\qquad + p(R+1)\,S\left(\sqrt{\frac{dp(R+1)}{p(R)+p(R+1)}}\,\right).
\end{aligned}
\end{equation*}
Recalling \eqref{IEdef}, \eqref{Sdef} and \eqref{Lubarv}, and rearranging terms, we find the expression in \eqref{IRid}.
\end{proof}

Note that $I^\dagger_R$ has a unique zero at $p$ given by
\begin{equation*}
p(0) = \tfrac12, \qquad p(k) = \tfrac12 (d-1)d^{-k}, \quad k = 1,\ldots,R, \qquad p(R+1) = \tfrac12d^{-R}. 
\end{equation*}
This shows that the fraction of the local time typically spent a distance $k$ away from the backbone decays exponentially fast in $k$.

%%%%%%%%%%% Appendix E %%%%%%%%%%%%%%%%%%%%%%%%%%

\section{Analysis of the upper variational formula}
\label{appE}

In this appendix we carry out the proof of the claims in Section~\ref{ss.varformup}, namely, we settle \eqref{IElberror} in Appendix \ref{IEid} and \eqref{limvarformlb} in Appendix \ref{varformlim}. The computations carried out in Appendix~\ref{appD} guide us along the way. 

%%%

\subsection{Identification of the rate function for the local times on the truncated tree}
\label{IEid}

To identify the rate function $I^\dagger_{E(\cW_R)}$ in Lemma~\ref{lem:VF}, we need to work out the two infima between braces in \eqref{VF}. The computation follows the same line of argument as in Appendix~\ref{appD}, but is more delicate. We will only end up with a lower bound. However, this is sufficient for the upper variational formula. 

To simplify the notation we write (recall Fig.~\ref{fig:foldedunit}):

\vspace{0.2cm}
\begin{tabular}{lll}
$(V_R,E_R)$ &=&  vertex and edge set of  $\cW_R$ \emph{without the tadpoles},\\
$\cO$ &=& \emph{top} vertex of $V_R$,\\
$\star$ &=& \emph{right-most bottom} vertex of $V_R$,\\ 
$\partial V_R$ &=& set of vertices at the \emph{bottom} of $V_R$,\\
$\Box$ &=& set of \emph{tadpoles},\\
$\Box_x$ &=& tadpole attached to $x \in \partial V_R\backslash\star$.
\end{tabular}

\vspace{0.2cm}\noindent
Note that $\partial V_R$ consists of $\star$ and the vertices to which the tadpoles are attached. Note that $\mathrm{int}(V_R) = V_R \setminus \partial V_R$ \emph{includes} $\cO$. 

\medskip\noindent
{\bf 1.}
Inserting \eqref{explaw} in Appendix~\ref{appB} into \eqref{rf1}--\eqref{rf2}, we get
\begin{equation*}
I^\dagger_{E(\cW_R)}(p) = (d+1) \sum_{x\in V_R} p(x)
+ \inf_{\beta \in (0,\infty)} \inf_{q \in \cP(V_R)} \sup_{\widehat{q} \in \cP(V_R)} L(\beta,q,\widehat{q} \mid p)
\end{equation*}
with
\begin{equation*}
L(\beta,q,\widehat{q} \mid p) =  - A - B - C - D,
\end{equation*}
where
\begin{equation*}
\begin{aligned}
A &= \sum_{x \in \mathrm{int}(V_R)} \beta q(x)\left\{1+\log\left(\frac{\sum_{y \sim x}\widehat{q}(y)}{\widehat{q}(x)}
\frac{p(x)}{\beta q(x)}\right)\right\},\\
B &= \sum_{x \in \partial V_R\backslash \star} \beta q(x)\left\{1+\log\left(\frac{\widehat{q}(x^\uparrow) 
+ d \widehat{q}(\Box_x)}{\widehat{q}(x)}\frac{p(x)}{\beta q(x)}\right)\right\},\\
C &= \beta q(\star) \left\{1+\log\left(\frac{\widehat{q}(\star^\uparrow) + d \widehat{q}(\cO)}{\widehat{q}(\star)}
\frac{p(\star)}{\beta q(\star)}\right)\right\},\\
D &= \sum_{x \in \Box}\beta q(x)\left\{\log\left(\frac{\widehat{q}(x^\uparrow)}{\widehat{q}(x)} \right)
- (\cL\lambda)\left(\frac{p(x)}{\beta q(x)}\right) \right\},
\end{aligned}
\end{equation*}
with $\cL\lambda$ the Legende transform of the cumulant generating function of $\psi$ (recall \eqref{cumulant}) and $x^\uparrow$ the unique vertex to which $x$ is attached upwards. (Recall that $y \sim x$ means that $x$ and $y$ are connected by an edge in $E_R$.) Note that $A,B,C$ each combine two terms, and that $A,B,C,D$ depend on $p$. We suppress this dependence because $p$ is fixed. 

\medskip\noindent
{\bf 2.}
Inserting the parametrisation $\widehat{q} = u/\|u\|_1$ and $q = v/\|v\|_1$ with $u,v\colon V_R \to (0,\infty)$ and putting $\beta q = v$, we may write  
\begin{equation}
\label{Krepr}
I^\dagger_{E(\cW^R)}(p) = (d+1) \sum_{x\in V_R} p(x) + \inf_{v\colon V_R \to (0,\infty)} \sup_{u\colon V_R \to (0,\infty)} L(u,v)
\end{equation}
with 
\begin{equation*}
L(u,v) = - A - B - C - D,
\end{equation*}
where 
\begin{equation}
\label{ABCD}
\begin{aligned}
A &= \sum_{x \in \mathrm{int}(V_R)} v(x)\left\{1+\log\left(\frac{\sum_{y \sim x}u(y)}{u(x)}\frac{p(x)}{v(x)}\right)\right\},\\
B &= \sum_{x \in \partial V_R \backslash \star} v(x)\left\{1+\log\left(\frac{u(x^\uparrow) 
+ d u(\Box_x)}{u(x)}\frac{p(x)}{v(x)}\right)\right\},\\
C &= v(\star) \left\{1+\log\left(\frac{u(\star^\uparrow) + d u(\cO)}{u(\star)}\frac{p(\star)}{v(\star)}\right)\right\},\\
D &= \sum_{x \in \Box}v(x)\left\{\log\left(\frac{u(x^\uparrow)}{u(x)} \right) - (\cL\lambda)\left(\frac{p(x)}{v(x)}\right) \right\}.
\end{aligned}
\end{equation}
Our task is to carry out the supremum over $u$ and the infimum over $v$ in \eqref{Krepr}.

\medskip\noindent
{\bf 3.}
First, we compute the infimum over $v$ for fixed $u$. (Later we will make a judicious choice for $u$ to obtain a lower bound.) Abbreviate
\begin{equation}
\label{ABCu}
\begin{aligned}
A_u(x) &= \frac{\sum_{y \sim x}u(y)}{u(x)}\,p(x), &x \in \mathrm{int}(V_R),\\
B_u(x) &= \frac{u(x^\uparrow) + d u(\Box_x)}{u(x)}\,p(x), &x\in \partial V_R\backslash\star,\\
C_u(\star) &= \frac{u(\star^\uparrow) + d u(\cO)}{u(\star)}\,p(\star). 
\end{aligned}
\end{equation}

\smallskip\noindent
$\bullet$
For $z \in V_R$, the first derivatives of $L$ are
\begin{equation*}
\begin{aligned}
&z \in \mathrm{int}(V_R)\colon
&\frac{\partial L(u,v)}{\partial v(z)} = -\log\left(\frac{A_u(z)}{v(z)}\right),\\
&z \in \partial V_R\backslash\star\colon
&\frac{\partial L(u,v)}{\partial v(z)} = -\log\left(\frac{B_u(z)}{v(z)}\right),\\
&z = \star\colon
&\frac{\partial L(u,v)}{\partial v(z)} = -\log\left(\frac{C_u(z)}{v(z)}\right),
\end{aligned}
\end{equation*}
while the second derivatives of $L$ equal $1/v(z)>0$. Hence the infimum is uniquely taken at 
\begin{equation*}
\begin{aligned}
&x \in \mathrm{int}(V_R)\colon
&\bar v(x) = A_u(x),\\
&x \in V_R \backslash \star \colon
&\bar v(x) = B_u(x),\\
&x = \star \colon
&\bar v(x) = C_u(x).
\end{aligned}
\end{equation*}

\smallskip\noindent
$\bullet$ For $z \in \Box$, the computation is more delicate. Define (see \eqref{mudef} in Appendix~\ref{appB})
\begin{equation*}
\mu(\alpha) = \alpha (\cL\lambda)^{'}(\alpha) - (\cL\lambda)(\alpha).
\end{equation*}
The function $\mu$ has range $(-\infty,\log \sqrt{d}\,]$, with the maximal value uniquely taken at $\alpha=\infty$. Therefore there are two cases.

\medskip\noindent
$\blacktriangleright$ $u(x)/u(x^\uparrow) \leq \sqrt{d}$:
Abbreviate $\alpha_u(z) = p(z)/v(z)$. For $z \in \Box$, 
\begin{equation*}
\begin{aligned}
\frac{\partial L(u,v)}{\partial v(z)} 
&= \log\left(\frac{u(z)}{u(z^\uparrow)}\right) 
+ (\cL\lambda)\left(\frac{p(z)}{v(z)}\right) - \frac{p(z)}{v(z)} (\cL\lambda)^{'}\left(\frac{p(z)}{v(z)}\right)\\
&= \log\left(\frac{u(z)}{u(z^\uparrow)}\right) - \mu(\alpha_u(z)), \\
\frac{\partial^2 L(u,v)}{v(z)^2} 
&=\frac{p^2(z)}{v^3(z)} (\cL\lambda)^{''}\left(\frac{p(z)}{v(z)}\right) >0,
\end{aligned}
\end{equation*}
where we use that $\cL\lambda$, being a Legendre transform, is strictly convex. Hence the infimum is uniquely taken at 
\begin{equation*}
\bar v(x) = \frac{p(x)}{\alpha_u(x)}, \qquad x \in \Box,
\end{equation*} 
with $\alpha_u(x)$ solving the equation
\begin{equation*}
\log\left(\frac{u(x)}{u(x^\uparrow)}\right) 
= \mu(\alpha_u(x)), \qquad x \in \Box. 
\end{equation*}
Since $\mu'(\alpha) = \alpha(\cL\lambda)''(\alpha)$ and $\cL\lambda$ is strictly convex (see Fig.~\ref{Fig:lambdaalt} in Appendix \ref{appB}), $\mu$ is strictly increasing and therefore invertible. Consequently,
\begin{equation*}
\alpha_u(x) = \mu^{-1}\left(\log\left(\frac{u(x)}{u(x^\uparrow)}\right)\right), \qquad x \in \Box.
\end{equation*} 
Putting the above formulas together, we arrive at (recall \eqref{ABCu})
\begin{equation}
\label{Lrepr}
\begin{aligned}
L(u) &= \inf_{v\colon V_R \to (0,\infty)} L(u,v)\\ 
&= - \sum_{x \in \mathrm{int}(V_R)} A_u(x) \quad - \sum_{x\in \partial V_R\backslash \star} B_u(x)\quad - C_u(\star)\quad
+ \sum_{x \in \Box} D_u(x)
\end{aligned}
\end{equation}
with (recall \eqref{ABCD})
\begin{equation*}
\begin{aligned}
D_u(x) &= - \frac{p(x)}{\alpha_u(x)}\left[\log\left(\frac{u(x^\uparrow)}{u(x)}\right) - (\cL\lambda)(\alpha_u(x))\right]\\
&= \frac{p(x)}{\alpha_u(x)} \big[(\cL\lambda)(\alpha_u(x)) - \mu(\alpha_u(x))\big]\\ 
&= p(x)\, (\cL\lambda)^{'}(\alpha_u(x)) 
= p(x)\, \big((\cL\lambda)^{'} \circ \mu^{-1}\big)\left(\log\left(\frac{u(x)}{u(x^\uparrow)}\right)\right).
\end{aligned}
\end{equation*}
In \eqref{cLrel2} in Appendix~\ref{appB} we show that $(\cL\lambda)' \circ \mu^{-1} = \lambda^{-1}$. Moreover In \eqref{lambdainverse} in Appendix~\ref{appB} we show that $(\lambda^{-1} \circ \log) = S$ with
\begin{equation*}
S(\beta) = d+1 - \beta - \frac{d}{\beta}, \qquad \beta \in (0,\sqrt{d}\,].
\end{equation*} 
Since $S$ has domain $(0,\sqrt{d}\,]$, $D_u(x)$ is only defined when $u(x)/u(x^\uparrow) \leq \sqrt{d}$, in which case
\begin{equation}
\label{Du}
D_u(x) = p(x)\,S\left(\frac{u(x)}{u(x^\uparrow)}\right), \qquad x \in \Box.
\end{equation}

\medskip\noindent
$\blacktriangleright$ $u(x)/u(x^\uparrow) > \sqrt{d}$: In this case $\frac{\partial L(u,v)}{\partial v(z)} > 0$, the infimum is uniquely taken at $\bar{v}(x)=0$, and
\begin{equation*}
D_u(x) = p(x)\,(\sqrt{d}-1)^2 = p(x)\,S(\sqrt{d}), \qquad x \in \Box,
\end{equation*}  
where we use \eqref{Legendrelim}. Note that the right-hand side does not depend on $u$.

\medskip\noindent
{\bf 4.} 
Next, we compute the supremum over $u$. The first derivatives of $L$ are 
\begin{equation}
\label{lines}
\begin{aligned}
&z \in \mathrm{int}(V_R) \backslash \cO\colon
&&\frac{\partial L(u)}{\partial u(z)} 
= \frac{\sum_{y \sim z} u(y)}{u^2(z)}\,p(z) - \sum_{y \sim z} \frac{1}{u(y)}\,p(y),\\
&z = \cO\colon
&&\frac{\partial L(u)}{\partial u(\cO)} 
= \frac{\sum_{y \sim \cO} u(y)}{u(\cO)^2}\, p(\cO) -\sum_{y: y^\uparrow = \cO} \frac{1}{u(y)}p(y) 
- \frac{d}{u(\star)}\,p(\star),\\
&z = \star\colon
&&\frac{\partial L(u)}{\partial u(\star)} 
= -\frac{1}{u(\cO)}\,p(\cO) + \frac{u(\star^\uparrow) + du(\cO)}{u(\star)^2}\,p(\star), \\
&z \in \partial V_R \backslash \star\colon
&&\frac{\partial L(u)}{\partial u(z)} 
= -\frac{1}{u(z^\uparrow)}\,p(z^\uparrow) + \frac{u(z^\uparrow)+du(\Box_z)}{u(z)^2}\,p(z)\\
&
&&\qquad \qquad + \left[\frac{u(\Box_z)}{u(z)^2} - \frac{d}{u(\Box_z)}\right]p(\Box_z) 
1_{\big\{\frac{u(z)}{u(z^\uparrow)}\leq \sqrt{d}\big\}}, \\
&z \in \Box\colon
&&\frac{\partial L(u)}{\partial u(z)} 
= -\frac{d}{u(z^\uparrow)}\,p(z^\uparrow) 
+ \left[-\frac{1}{u(z^\uparrow)} +\frac{du(z^\uparrow)}{u(z)^2}\right]\,p(z)\,
1_{\big\{\frac{u(z)}{u(z^\uparrow)}\leq \sqrt{d}\big\}}.
\end{aligned}
\end{equation}
The second derivates of $L$ are all $<0$. The first line in \eqref{lines} can be rewritten as
\begin{equation*}
\sum_{y \sim z} u(y) \left[\frac{p(z)}{u^2(z)} - \frac{p(y)}{u^2(y)} \right],
\end{equation*}
which is zero when
\begin{equation}
\label{ubasic}
\bar u(x) = \sqrt{p(x)},  \qquad x \in V_R.
\end{equation}
Given the choice in \eqref{ubasic}, the fifth line in \eqref{lines} is zero when
\begin{equation}
\label{utadpole}
\bar u(x) = \sqrt{\frac{dp(x^\uparrow)p(x)}{dp(x^\uparrow)+p(x)}}, \qquad x \in \Box.
\end{equation} 
Indeed, the derivative is strictly negative when the indicator is $0$ and therefore the indicator must be $1$. But the latter is guaranteed by \eqref{ubasic}--\eqref{utadpole}, which imply that
\begin{equation*}
\frac{\bar u(x)}{\bar u(x^\uparrow)} = \sqrt{\frac{dp(x)}{dp(x^\uparrow)+p(x)}} \leq \sqrt{d}, \qquad x \in \Box. 
\end{equation*}
Given the choice in \eqref{ubasic}--\eqref{utadpole}, also the fourth line in \eqref{lines} is zero. Thus, only the second and third line in \eqref{lines} are non-zero, but this is harmless because $\cO,\star$ carry a negligible weight in the limit as $R \to \infty$ because of the constraint $p(\partial V_R \cup \cO) \leq 1/R$ in Lemma~\ref{lem:VF} (recall \eqref{parRdef}). 

Inserting \eqref{ubasic}--\eqref{utadpole} into \eqref{Lrepr} and using \eqref{ABCu}, \eqref{Du}, we get the following lower bound:
\begin{equation*}
\begin{aligned}
&\sup_{u\colon V_R \to (0,\infty)} L(u)\\
&\geq - \sum_{x \in \mathrm{int}(V_R)} A_{\bar u}(x) \quad 
- \sum_{x\in \partial V_R\backslash \star} B_{\bar u}(x)\quad 
- C_{\bar u}(\star) + \sum_{x \in \Box} D_{\bar u}(x)\\
&= - \sum_{x \in \mathrm{int}(V_R)} \sum_{y \sim x} \sqrt{p(y)p(x)} 
- \sum_{x\in \partial V_R \backslash \star} \sqrt{p(x)}\left(\sqrt{p(x^\uparrow)} 
+ d\sqrt{\frac{dp(x)p(\Box_x)}{dp(x)+p(\Box_x)}}\right) \\
&\quad -\sqrt{p(\star)}\left(\sqrt{p(\star^\uparrow)}+ d\sqrt{p(\cO)}\right)\\
&\quad + \sum_{x \in \Box} p(x)\,\left(d+1-\sqrt{d} \left[\sqrt{\frac{p(x)}{d p(x^\uparrow) + p(x)}}
+ \sqrt{\frac{d p(x^\uparrow) + p(x)}{p(x)}}\,\right]\right).
\end{aligned}
\end{equation*}

\medskip\noindent
{\bf 5.}
Using the relation $(d+1) p(x) = \sum_{y\sim x} p(x)$, $x\in \mathrm{int}(V_R)$, we get from \eqref{Krepr} that
\begin{equation*}
I^\dagger_{E(\cW^R)}(p) \geq K^1_R(p) + K^2_R(p) 
\end{equation*}
with
\begin{equation*}
\begin{aligned}
K^1_R(p) 
&= \sum_{x \in \mathrm{int}(V_R)} \sum_{y \sim x} \left[p(x) - \sqrt{p(x)p(y)}\,\right] \\
&= \sum_{\{x,y\} \in \widehat{E}_R} \left(\sqrt{p(x)} - \sqrt{p(y)}\,\right)^2 
+ \left[p(\cO)-\sqrt{p(\cO)p(\star)}\,\right] - \sum_{x\in \partial V_R}\left[ p(x) - \sqrt{p(x)p(x^\uparrow)}\,\right]
\end{aligned}
\end{equation*}
and
\begin{equation*}
\begin{aligned}
K^2_R(p) 
&= \sum_{x\in \partial V_R \backslash \star}\left[(d+1) p(x) - \sqrt{p(x)}\left(\sqrt{p(x^\uparrow)} 
+ d\sqrt{\frac{dp(x)p(\Box_x)}{dp(x)+p(\Box_x)}}\right)\right]\\
&\qquad + (d+1) p(\star)-\sqrt{p(\star)}\left(\sqrt{p(\star^\uparrow)} + d\sqrt{p(\cO)}\right) \\
&\qquad + \sum_{x \in \Box} p(x)\,\left[d+1-\sqrt{d}\,\left(\sqrt{\frac{p(x)}{d p(x^\uparrow) + p(x)}}
+ \sqrt{\frac{d p(x^\uparrow) + p(x)}{p(x)}}\,\right)\right].
\end{aligned}
\end{equation*}
The first sum in the right-hand side of $K^1_R(p)$ equals the \emph{standard rate function} $I_{\widehat{E}_R}(p)$ given by \eqref{e:defIJ}, with 
\begin{equation*}
\widehat{E}_R = E_R \setminus \{\cO,\star\}
\end{equation*} 
the set of edges in the unit $\cW_R$ \emph{without the tadpoles} and \emph{without the edge} $\{\cO,\star\}$ (i.e., $\widehat{E}_R = E(\cT^*_R)$; recall Fig.~\ref{fig:unit}). Rearranging and simplifying terms, we arrive at
\begin{equation}
\label{lbfinal1}
I^\dagger_{E(\cW^R)}(p) \geq  I_{\widehat{E}_R}(p)+ K^3_R(p) 
\end{equation}
with
\begin{equation*}
K^3_R(p) = S_{\partial V_R \backslash \star}(p) + S_{\cO,\star}(p) + S_{(\partial V_R \backslash \star) \cup \Box}(p),
\end{equation*}
where
\begin{equation}
\label{lbfinal3}
\begin{aligned}
S_{\partial V_R \backslash \star}(p) 
&= d \sum_{x\in \partial V_R \backslash \star} p(x),\\
S_{\cO,\star}(p) 
&= \left(\sqrt{p(\cO)} - \sqrt{p(\star)}\right)^2 + (d-1)\big[p(\star) - \sqrt{p(\cO)p(\star)}\,\big],\\
S_{(\partial V_R \backslash \star) \cup \Box}(p) 
&= - \sum_{x\in \partial V_R \backslash \star} p(x)\,d\sqrt{\frac{dp(\Box_x)}{dp(x)+p(\Box_x)}}\\
&\qquad +  \sum_{x\in \partial V_R \backslash \star} p(\Box_x)\,\left(d+1-\sqrt{d}\,\left[\sqrt{\frac{p(\Box_x)}{d p(x) + p(\Box_x)}}
+ \sqrt{\frac{d p(x) + p(\Box_x)}{p(\Box_x)}}\,\right]\right).
\end{aligned}
\end{equation}

\medskip\noindent
{\bf 6.}
Since $\sqrt{p(\cO)p(\star)} \leq \tfrac12[p(\cO)+p(\star)]$, the boundary constraint $\sum_{x\in \partial V_R \cup \cO} p(x) \leq 1/R$ implies that $S_{\partial V_R \backslash \star}(p) + S_{\cO,\star}(p) = O(1/R)$. The same constraint implies that the first sum in $S_{(\partial V_R \backslash \star) \cup \Box}(p) $ is $O(1/R)$. Hence 
\begin{equation*}
K^3_R(p) = O(1/R) + \sum_{x\in \partial V_R \backslash \star} p(x)\,F\left(\tfrac{p(\Box_x)}{p(x)}\right)
\end{equation*} 
with 
\begin{equation*}
F(w) = w \left(d+1-\sqrt{d}\,\left[\sqrt{\frac{w}{d+w}} + \sqrt{\frac{d+w}{w}}\,\right]\right).
\end{equation*}
The map $w \mapsto F(w)$ is continuous on $(0,\infty)$ with 
\begin{equation*}
F(w) = \left\{\begin{array}{lll}
&\sqrt{w} + (d+1)w + O(w^{3/2}), &w \downarrow 0,\\
&[(d+1)-2\sqrt{d}\,]\, w + O(w^{-1}), &w \to \infty.
\end{array}
\right.
\end{equation*}
From this we see that if $d \geq 4$, then there exists a $C \in (1,\infty)$ such that
\begin{equation}
\label{Fbd}
F(w)+C \geq \left(1-\sqrt{w}\,\right)^2, \qquad w \in [0,\infty).
\end{equation} 
Hence we have the lower bound
\begin{equation*}
\begin{aligned}
K^3_R(p) 
&\geq O(1/R) + \sum_{x\in \partial V_R \backslash \star}
p(x) \left[-C + \left(1-\sqrt{\tfrac{p(\Box_x)}{p(x)}}\,\right)^2\right]\\
&= O(1/R) + \sum_{x\in \partial V_R \backslash \star} \left(\sqrt{p(x)}-\sqrt{p(\Box_x)}\,\right)^2.
\end{aligned}
\end{equation*}
Via \eqref{lbfinal1}--\eqref{lbfinal3}, it follows that
\begin{equation}
\label{IEbd}
I^\dagger_{E(\cW^R)}(p) \geq O(1/R) + I_{\widetilde{E}_R}(p), \qquad R \in \N,
\end{equation}
with $I_{\widetilde{E}_R}(p)$ the \emph{standard rate function} given by \eqref{e:defIJ}, with 
\begin{equation*}
\widetilde{E}_R = \widehat{E}_R \cup \big[\cup_{x \in \partial V_R \setminus \star} \{x,\Box_x\}\big]
\end{equation*} 
the set of edges in the unit $\widetilde\cW_R$ that is obtained from the unit $\cW_R$ by \emph{removing the edge} $\{\cO,\star\}$ (i.e., $\widetilde{E}_R = E(\widetilde\cW_R)$; recall Fig.~\ref{fig:foldedunit}). This completes the proof of \eqref{IElberror}.
 
\begin{remark}
\label{bdcontrol}
{\rm The condition $d \geq 4$ is needed only in \eqref{Fbd}. For $d=2,3$ we have $F(w)+C \geq \theta_c(1-\sqrt{w}\,)^2$ with $\theta_c = d+1-2\sqrt{d} \in (0,1)$. Consequently, the edges $\{x,\Box_x\}$, $x \in \partial V_R\setminus \star$, carry a weight that is smaller than that of the edges in $\cT$, which may cause the optimal $p$ to stick to the boundary as $R\to\infty$, in which case we do not have \eqref{IEbd}.} \hfill$\spadesuit$
\end{remark} 
 
%%%

\subsection{Limit of the upper variational formula}
\label{varformlim}

Note that 
\begin{equation*}
\widetilde\cW_R \subseteq \cT,
\end{equation*}
with $\cT$ the infinite tree. Consequently,  
\begin{equation*}
I_{\widetilde{E}_R}(p) = I_{E(\cT)}(p) - (d-1) \sum_{x \in \partial V_R \setminus \star} p(x),   
\qquad \forall\,p \in \cP(V(\cT))\colon\,\supp(p) = V(\widetilde\cW_R), 
\end{equation*}
where the sum compensates for the contribution coming from the edges in $\cT$ that link the vertices in $\partial V_R \setminus \star$ to the vertices one layer deeper in $\cT$ that are not tadpoles. Since this sum is $O(1/R)$, we obtain (recall \eqref{VF})
\begin{equation*}
\begin{aligned}
\chi^+_R(\varrho) &= \inf_{ {p \in \cP(V(\cW_R))\colon} \atop {p(\partial_{\cW_R}) \leq 1/R}}
\left\{I^\dagger_{E(\cW_R)}(p) + \varrho J_{V(\cW_R)}(p)\right\}\\
&\geq O(1/R) + \inf_{ {p \in \cP(V(\cT))\colon} \atop {\supp(p) = V(\widetilde\cW_R), 
\,p(\partial_{\widetilde\cW_R}) \leq 1/R}}
\left\{I_{E(\cT)}(p) + \varrho J_{V(\cT)}(p)\right\}\\
&\geq O(1/R) + \chi_{\cT}(\rho),
\end{aligned}
\end{equation*}
where the last inequality follows after dropping the constraint under the infimum and recalling \eqref{varform}. This completes the proof of \eqref{limvarformlb}.

%%%%%%%%%%%% REFERENCES %%%%%%%%%%%%%%%%%%%%%%%%%%%%%%%%%%%%

%%%%%%%%%%%%%%%%%%%%%%%%%%%%%%%%%%%%%%%%%%%%%%

\end{document}